\documentclass[ onefignum,onetabnum]{siamart190516}

\usepackage{times}
\usepackage{microtype}
\usepackage{graphicx}
\usepackage{booktabs} %
\usepackage{subfigure}
\usepackage{hyperref}

\usepackage{color}
\usepackage{url}            %
\usepackage{amsfonts}       %
\usepackage{nicefrac}       %
\usepackage{microtype}      %

\usepackage{multirow}

\usepackage{times}
\usepackage{algorithm}
\usepackage{xspace}

\usepackage{amsfonts}
\usepackage{graphicx}
\usepackage{amsmath}
\usepackage{amssymb}
\usepackage{bbding}
\usepackage{pifont}
\usepackage{wasysym}
\usepackage{bm}
\usepackage{algorithmic}
\usepackage{epstopdf}
\usepackage{psfrag}
\usepackage{bbm}
\usepackage{enumerate}
\usepackage{subeqnarray}
\usepackage{cases}

\def\circlenum#1{\text{\textcircled{\textbf{#1}}}}

\newcommand{\argmin}{\mathop{\rm arg\min}}

\newcommand{\bbR}{\mathbb{R}}

\newcommand{\Range}{{\mathcal{R}}}

\renewcommand{\d}{{\rm{d}}}

\usepackage{color}
\def\rc{\color{red}}

\usepackage{soul}

\newcommand{\beq}{\begin{equation}}
\newcommand{\eeq}{\end{equation}}
\newcommand{\beqa}{\begin{eqnarray}}
\newcommand{\eeqa}{\end{eqnarray}}
\newcommand{\beqas}{\begin{eqnarray*}}
\newcommand{\eeqas}{\end{eqnarray*}}
\newcommand{\bi}{\begin{itemize}}
\newcommand{\ei}{\end{itemize}}
\newcommand{\ba}{\begin{array}}
\newcommand{\ea}{\end{array}}

\def\argmin{{\rm argmin}}

\def\exp{{\rm exp}}

\def\cA{{\cal A}}

\usepackage{lipsum}
\usepackage{amsfonts}
\usepackage{graphicx}
\usepackage{epstopdf}
\usepackage{algorithmic}
\ifpdf
  \DeclareGraphicsExtensions{.eps,.pdf,.png,.jpg}
\else
  \DeclareGraphicsExtensions{.eps}
\fi

\newsiamremark{remark}{Remark}
\newsiamremark{example}{Example}
\newsiamremark{hypothesis}{Hypothesis}
\newsiamremark{assumption}{Assumption}
\crefname{hypothesis}{Hypothesis}{Hypotheses}
\newsiamthm{claim}{Claim}

\headers{Unified Acceleration of High-order Algorithms}{Chaobing Song, Yong Jiang, and Yi Ma}

\title{Unified Acceleration of  High-Order Algorithms under \\ H\"{o}lder Continuity and Uniform Convexity\thanks{Submitted\funding{This work was done during Chaobing Song's visit to Professor Yi Ma's group at UC Berkeley. The work is partially supported by the TBSI program and EECS Startup fund of Professor Yi Ma.}}
}

\author{Chaobing Song\thanks{Tsinghua-Berkeley Shenzhen Institute (TBSI), Tsinghua University
  (\email{songcb16@mails.tsinghua.edu.cn}).}
  \and Yong Jiang\thanks{Tsinghua-Berkeley Shenzhen Institute (TBSI), Tsinghua University (\email{jiangy@sz.tsinghua.edu.cn})}
\and Yi Ma\thanks{EECS Department, University of California, Berkeley
  (\email{yima@eecs.berkeley.edu}).}
}

\usepackage{amsopn}

\ifpdf
\hypersetup{
  pdftitle={Towards Unified Acceleration of  High-Order Algorithms under \\ H\"{o}lder Continuity and Uniform Convexity},
  pdfauthor={Chaobing Song, Yong Jiang and Yi Ma}
}
\fi

\begin{document}

\maketitle

\begin{abstract}
{ In this paper, through a very intuitive {\em vanilla proximal method} perspective,  we derive accelerated high-order optimization algorithms for minimizing a convex function that has
H\"{o}lder continuous derivatives. In this general convex setting, we propose a concise {\em unified acceleration framework (UAF)}, which reconciles the two different high-order acceleration approaches, one by Nesterov and Baes \cite{nesterov2008accelerating, baes2009estimate,nesterov2018implementable} and one by Monteiro and Svaiter \cite{monteiro2013accelerated}. As result, the UAF unifies the high-order acceleration instances \cite{
nesterov2008accelerating,baes2009estimate, nesterov2018implementable,grapiglia2019accelerated,
grapiglia2019tensor, monteiro2013accelerated,jiang2018optimal, bubeck2018near,gasnikov2018global} of the two approaches by only two problem-related parameters and  two additional parameters for framework design. 
 Furthermore, the UAF (and its analysis) is \emph{the first approach} to make high-order methods applicable for high-order smoothness conditions with respect to non-Euclidean norms. 
If the function is further uniformly convex, we propose a {\em general restart scheme} for the UAF. The iteration complexities of instances of both the UAF and the restarted UAF  match existing lower bounds in most important cases \cite{arjevani2017oracle,grapiglia2019tensor}.   
 For practical implementation, we introduce a new and effective heuristic that significantly simplifies the binary search procedure required by the framework. We use experiments to verify the effectiveness of the heuristic and demonstrate clear and consistent advantages of high-order acceleration methods over first-order ones, in terms of run-time complexity.
 Finally, the UAF is proposed directly in 
 the general composite convex setting, thus show that the existing high-order algorithms \cite{
nesterov2008accelerating,baes2009estimate, nesterov2018implementable,
grapiglia2019tensor, bubeck2018near,gasnikov2018global} can be naturally extended to the general composite convex setting.
 }
\end{abstract}



\section{Introduction}\label{sec:intro}
In optimization, people often consider the problem of minimizing a convex function: 
\begin{equation}
\min_{x\in\bbR^d} f(x).\label{eq:smooth-prob}
\end{equation}
A typical assumption is that $f(x)$ has $L$-Lipschitz continuous gradients with respect to ($w.r.t.$) the Euclidean norm $\|\cdot\|_2$, $i.e.,$
\begin{eqnarray}
\|\nabla f(x) - \nabla f(y)\|_2\le L \|x - y\|_2,\label{eq:lip}
\end{eqnarray}
where $L>0$ is the Lipschitz constant. For this problem, to find an $\epsilon$-accurate solution $x$ such that $f(x) - f(x^*)\le \epsilon$, the classic gradient descent method:
$$
x_{k+1} = x_k - \eta \nabla f(x_k)
$$
with $\eta \le 1/L$ takes $O(\epsilon^{-1})$ iterations. Nevertheless, it is known that from \cite{nesterov1998introductory}, for a convex function $f(x)$ with $L$-Lipschitz continuous gradients, a lower-bound for the number of iterations for any first-order algorithms is known to be 
\begin{equation}
O\big(\epsilon^{-1/2}\big), \quad \text{($L$-Lipschitz continuous gradients)}.
\label{bound:agd}
\end{equation}In the seminal work \cite{nesterov1983method}, Nesterov has introduced an acceleration technique, the so-called \emph{accelerated gradient descent (AGD)} algorithm, that achieves this optimal lower bound. This algorithm dramatically improves the convergence rate of smooth convex optimization with negligible per-iteration cost. { Besides the smooth convex problem \eqref{eq:smooth-prob} under the Euclidean norm setting \eqref{eq:lip}, AGD can also be generalized to solve the composite convex problem \cite{beck2009fast,allen2017katyusha,diakonikolas2019approximate}, in which the objective function may contain a second possibly non-smooth but simple convex term (see \eqref{eq:prob}). Meanwhile, AGD can be extended to the more general (non-Euclidean) norm settings \cite{diakonikolas2017accelerated,allen2017katyusha,krichene2015accelerated}, also achieving the optimal rate \eqref{bound:agd}.}

\subsection{High-order Acceleration Methods with Lipschitz Continuity}\label{sec:Lip}
To hope for a better iteration complexity beyond $O(\epsilon^{-1/2})$,  $f(x)$ needs to be smooth for its high-order derivatives. A common assumption is that $f(x)$ has $(p, \nu, L)$-H\"{o}lder continuous derivatives: 
\begin{eqnarray}
\frac{1}{(p-1)!}\|\nabla^p f(x) - \nabla^p f(y)\|_{2}\le L \|x - y\|_2^{\nu},\label{eq:deri-lip-holder}
\end{eqnarray}
for some $\nu \in [0, 1], p\in\mathbb{Z}_+$. Notice that for $p = 1$ and $\nu = 1$, this condition becomes the first order $L$-Lipschitz continuous gradient \eqref{eq:lip} above.  Here, for $p \ge 2$, the $\|\cdot\|_{2}$ norm of a $p$-th order tensor denotes its operator norm \cite{nesterov2018implementable} $w.r.t.$ the vector 2-norm $\|\cdot\|_2$. Sometimes, when $\nu = 1$, the function is said to have $(p, L)$-Lipschitz continuous derivatives:
\begin{eqnarray}
\frac{1}{(p-1)!}\|\nabla^p f(x) - \nabla^p f(y)\|_2\le L \|x - y\|_2.\label{eq:lip-holder}
\end{eqnarray}
In general, if we were able to utilize higher-order derivatives with $p \ge 2$, we expect to obtain algorithms with higher convergence rates. The higher is $p$ (and $\nu$), the higher the rate could be.

If a convex function $f(x)$ has $(p, L)$-Lipschitz continuous derivatives \eqref{eq:lip-holder}, the recent work \cite{arjevani2017oracle} has given a lower-bound on the complexity: any deterministic algorithm would need at least \begin{equation}
O\big(\epsilon^{-\frac{2}{3p+1}}\big), \quad \text{($(p, L)$-Lipschitz continuous derivatives)}
\label{bound:pL-lip}
\end{equation} iterations to find an $\epsilon$-accurate solution. For the special case $p=2$, \cite{nesterov2008accelerating} has proposed an ``{\em accelerated cubic regularized Newton method}'' (ACNM) that achieves an iteration complexity of $O(\epsilon^{-\frac{1}{3}})$. From a  different approach, after proposing an accelerated hybrid proximal extragradient (A-HPE) framework, \cite{monteiro2013accelerated} has implemented an ``{\em accelerated Newton proximal extragradient}'' (A-NPE) instance of the A-HPE framework that has achieved the optimal complexity  $O(\epsilon^{-\frac{2}{7}})$\footnote{ When talking about iteration complexity, we mean the complexity in terms of outer iteration without concerning about the inner implementation of subproblems.} for $p = 2$, although each iteration requires a nontrivial binary search procedure. 

To achieve better complexity results and also being encouraged by the fact that third-order methods can often be implemented as efficiently as second-order methods \cite{nesterov2018implementable}, there is an increasing interest to extend ACNM  and implement the A-HPE framework to even higher-order smoothness settings $(p\in\{3,4\ldots,\})$ \cite{baes2009estimate, nesterov2018implementable,jiang2018optimal,gasnikov2018global,bubeck2018near}. In particular, by extending ACNM, \cite{baes2009estimate} and \cite{nesterov2018implementable} have proposed accelerated tensor methods with $O(\epsilon^{-\frac{1}{p+1}})$ iteration complexity for $p\in\{2,3,\ldots\}$. Meanwhile, by implementing A-HPE, \cite{monteiro2013accelerated, jiang2018optimal,gasnikov2018global,bubeck2018near} have proposed accelerated methods that achieve the optimal $O(\epsilon^{-\frac{2}{3p+1}})$ iteration complexity, although just like A-NPE, all these methods need the nontrivial binary search procedure. 

Hence the current situation seems to be: methods \cite{ nesterov2008accelerating,baes2009estimate,nesterov2018implementable,grapiglia2019accelerated,grapiglia2019tensor} by extending ACNM  have advantages with simpler implementation, while methods \cite{monteiro2013accelerated,jiang2018optimal,gasnikov2018global,bubeck2018near} by implementing  A-HPE
 can in theory achieve the optimal rate $O(\epsilon^{-\frac{2}{3p+1}})$. However, it remains somewhat mysterious how we could reconcile the differences between these two approaches. { In addition, the A-HPE framework is somewhat abstract so implementing it in the high-order setting requires rather nontrivial techniques \cite{monteiro2013accelerated, jiang2018optimal,gasnikov2018global,bubeck2018near}. It remains unclear how to propose a concise but equivalently powerful alternative to the A-HPE framework and obtain these different instances of A-HPE in a unified way. Furthermore, although AGD can be generalized to general non-Euclidean norm settings, up to now, it is not known whether high-order methods can have a similar generalization. Finally, both the ACNM and A-HPE approaches do not directly address the composite convex setting (see \eqref{eq:prob}) at the framework level, hence obtaining high-order algorithms in this setting is highly desired and seems nontrivial  \cite{monteiro2013accelerated,grapiglia2019accelerated,jiang2018optimal}.

}
\subsection{Acceleration under H\"{o}lder Continuity and Our Results}
Besides the Lipschitz continuous setting, the more general H\"{o}lder continuous setting \eqref{eq:deri-lip-holder} is also of increased interest, partly for designing universal optimization schemes \cite{nesterov2015universal,yurtsever2015universal,grapiglia2019tensor,cartis2019universal}. If $f(x)$ has $(1, \nu, L)$-H\"{o}lder continuous gradients, a lower bound for the iteration complexity is known to be \cite{nemirovsky1983problem}: 
\begin{equation}
   O\Big(\epsilon^{-\frac{2}{1+3\nu}}\Big), \quad \text{($(1, \nu, L)$-H\"{o}lder continuous gradients).} 
  \label{bound:hol-grad}
\end{equation}  An algorithm that can achieve this lower bound has been proposed in \cite{nemirovskii1985optimal}. 

For the more general setting of $(p, \nu, L)$-H\"{o}lder continuous derivatives, during the preparation of this paper,  \cite{grapiglia2019tensor} has given a lower bound of iteration complexity \begin{equation}
    O\Big(\epsilon^{-\frac{2}{3(p+\nu)-2}}\Big), \quad \text{($(p, \nu, L)$-H\"{o}lder continuous derivatives).}
    \label{bound:hol-der}
\end{equation} By extending Nesterov's method in \cite{nesterov2018implementable}, \cite{grapiglia2019tensor} has proposed a method that achieves the iteration complexity $O(\epsilon^{-\frac{1}{p+\nu}})$. To the best of our knowledge, methods that can achieve the lower bound $O\big(\epsilon^{-\frac{2}{3(p+\nu)-2}}\big)$ are still unknown.

In this paper, for the minimization of convex functions with $(p, \nu, L)$-H\"{o}lder continuous derivatives, we propose a {\em unified acceleration framework} (UAF), see Algorithm \ref{alg:dis}, that achieves the  iteration complexity of $O\big(\epsilon^{-\frac{2}{3(p+\nu)-2}}\big)\;(p\in\{1,2,\ldots\},  \nu\in[0,1]$ with $p+\nu\ge 2$, $L>0$), which matches the lower bound \cite{grapiglia2019tensor}. To be more precise, if a convex function $f(x)$ has $(p, \nu, L)$-H\"{o}lder continuous derivatives, our algorithm can find an $\epsilon$-accurate solution with
\begin{eqnarray}
O\left(\epsilon^{-\frac{q}{(q+1)(p+\nu)-q}}\right) \label{eq:rate}
\end{eqnarray}
iterations,
where $q$ is a tunable parameter\footnote{As we will later see, $q$ is the order of the uniform convexity of the proxy-function for framework design. \cite{grapiglia2019tensor} has used a uniformly convex proxy-function with $q=(p+\nu)$-th order, while  \cite{jiang2018optimal,gasnikov2018global,bubeck2018near} have used a uniformly convex proxy-function with $q=2$-nd order.} such that $2\le q\le p+\nu$. Notice that our result and algorithm unify previously known results as (important) special cases: 
\begin{itemize}
\item For the case of $L$-Lipschitz continuous gradients \cite{nesterov1998introductory} where  $p=\nu=1$ and $q=2$, the rate \eqref{eq:rate} of the proposed algorithm achieves the lower bound $O\big(\epsilon^{-\frac{1}{2}}\big)$ of \eqref{bound:agd}. 
\item For the more general setting of $(p, \nu, L)$-H\"{o}lder continuous derivatives: when $p\in\{2,3,\ldots,\}, q = p+\nu$, it recovers the complexity  $O\big(\epsilon^{-\frac{1}{p+\nu}}\big)$ of the method  in \cite{grapiglia2019tensor}. 
\end{itemize}
Meanwhile, by setting $q=2$,  the rate \eqref{eq:rate} of the UAF is \emph{the first} convergence result that matches the lower bound $O\big(\epsilon^{-\frac{2}{3(p+\nu)-2}}\big)$ of \eqref{bound:hol-der} \cite{grapiglia2019tensor} under the general H\"{o}lder continuous setting.

{
Besides the unified convergence rate \eqref{eq:rate}, the UAF has several significant improvements over the ACNM approach and the A-HPE framework.  %
First, the UAF provides a continuous transition from the ACNM approach to the A-HPE framework by choosing $q$ from $p+\nu$ to $2$. 
Second, as we will soon see, the UAF can be conveniently instantiated by only specifying
two problem-related parameters and two adjustable parameters for framework design, and thus recover the high-order acceleration algorithms  \cite{
nesterov2008accelerating,baes2009estimate, nesterov2018implementable,grapiglia2019accelerated,
grapiglia2019tensor, monteiro2013accelerated,jiang2018optimal, bubeck2018near,gasnikov2018global} without extra effort.
Third, we provide \emph{the first} and also a unified convergence rate analysis for both the Euclidean and non-Euclidean norm settings, and thus opens the possibility of applying high-order methods in the non-Euclidean norm setting.\footnote{which is pertinent to many important practical problems such as logistic regression loss in machine learning, see Example \ref{exam:smooth}.} Fourth, the UAF is proposed and analyzed directly under the  composite convex setting (see \eqref{eq:prob}), hence our results imply that all existing high-order algorithms \cite{
nesterov2008accelerating,baes2009estimate, nesterov2018implementable,
grapiglia2019tensor, bubeck2018near,gasnikov2018global} can be naturally extended to the general composite convex setting.

Last but not the least, 
in the high-order setting, to obtain the optimal rate that matches the lower bound \cite{arjevani2017oracle}, 
we must employ a binary search procedure to find a suitable coupling coefficient in each iteration, which may substantially slow down the practical performance \cite{nesterov2018implementable}. Therefore, in addition to the above theoretical results, we introduce a simple heuristic for finding the coupling coefficient, suggested by our analysis, so that the resulting implementation does not need a binary search procedure required by the optimal acceleration method. Our experiments show that this simple heuristic is extremely effective and can easily ensure the conditions needed to achieve the optimal rate. This leads to a very practical implementation of the optimal acceleration algorithms without extra implementation cost, alleviating concerns raised by \cite{nesterov2018implementable}.

}

\subsection{Acceleration under Uniform Convexity and Our Results}
The above result for optimal complexity \eqref{eq:rate} is given for the general convexity setting with $(p,\nu, L)$-H\"{o}lder continuous derivatives. When $f(x)$ has additional nice properties such as {\em uniform convexity}, we should expect even better complexity. To be more precise, assume that $f(x)$ is $(s,\sigma)$-uniformly convex, {\em i.e.}:
\begin{eqnarray}
f(x)\ge f(y)+\langle \nabla f(y), x-y\rangle + \frac{\sigma}{s}\|x-y\|_2^{s},
\end{eqnarray}
for $s\ge 2, \sigma>0$, where $(2,\sigma)$-uniform convexity is also known as $\sigma$-strong convexity. It is known that for $s=2, \sigma >0$,
when $f(x)$ having $L$-Lipschitz continuous gradient, \cite{nesterov1998introductory} has provided a lower bound for the iteration complexity 
\begin{equation}
O\left(\sqrt{\frac{L}{\sigma}}\log\frac{1}{\epsilon}\right), \quad \text{($L$-Lipschitz continuous gradient,  $\sigma$-strong convexity)}.
\end{equation}
When $f(x)$ has $(2, L)$-Lipschitz continuous derivatives, \cite{arjevani2017oracle} has provided a lower bound for the iteration complexity \begin{equation}
    O\Big(\Big(\frac{L}{\sigma}\Big)^{\frac{2}{7}}+\log\log\Big( \frac{\sigma^3}{L^2\epsilon}\Big)\Big),\; \text{($(2,L)$-Lipschitz {continuous derivatives},  $\sigma$-strong convexity)},
    \label{bound:lip-uc}
\end{equation} and it has also proposed a method based on restarting A-HPE \cite{monteiro2013accelerated} that achieves a complexity upper-bounded by $O\left(\left(\frac{L}{\sigma}\right)^{2/7}\log\frac{L^2}{\sigma^3} + \log\log\left(\frac{\sigma^3}{L^2\epsilon}\right)\right)$, quite close to the lower bound.

In this paper, we show that in the uniformly convex setting, the idea of restart for ACNM in \cite{nesterov2008accelerating} is also applicable to our algorithm and can significantly improve the iteration complexity \eqref{eq:rate}. Inspired by that work,  in this paper, we introduce a more {\em general restart scheme}, see Algorithm \ref{alg:restart}, that is applicable for accelerating almost all uniformly convex optimization algorithms.

We show that for $(s, \sigma)$-uniformly convex and $(p, \nu, L)$-H\"{o}lder continuous functions,
if $s = p+\nu$, then the UAF  with the proposed restart scheme applied, needs at most
\begin{eqnarray}
O\left(\left( \frac{L}{\sigma}\right)^{\frac{q}{(q+1)(p+\nu)-q}}\log \frac{1}{\epsilon}\right) \label{eq:linear-res}
\end{eqnarray} 
iterations to find an $\epsilon$-accurate solution, where $q \in [2, p+\nu]$ is the tunable parameter as before.
If $s<p+\nu$, then the resulting algorithm needs at most
\begin{eqnarray}
O\left(\left(\frac{L}{\sigma} \right)^{\frac{q}{(q+1)(p+\nu)-q}}+\log\log\left(\left(\frac{\sigma^{p+\nu}}{L^s} \right)^{\frac{1}{p+\nu-s}}\frac{1}{\epsilon}\right)\right)  \label{eq:super-res}
\end{eqnarray}
iterations. If $s>p+\nu$, then the algorithm needs at most 
\begin{eqnarray}
O\left(\left(  
\frac{L}{\sigma}\right)^{\frac{q}{(q+1)(p+\nu)-q}}\left(\frac{\sigma}{\epsilon}\right)^{\frac{(s-p-\nu)q}{s((q+1)(p+\nu)-q)}}\right)
\label{eq:sublinear-res}
\end{eqnarray}
iterations.

Notice that according to \eqref{eq:linear-res}, when $p=\nu=1, s=2$, with the design parameter $q=2$, we recover the optimal rate of accelerated gradient descent (AGD) in the strong convex setting \cite{nesterov1998introductory}. According to \eqref{eq:super-res}, when $p=2, \nu=1, s = 2$, with $q=2$, our algorithm eliminates the logarithmic factor in the first term of the upper bound given in \cite{arjevani2017oracle} and achieves the iteration complexity  $O\big(\left(\frac{L}{\sigma}\right)^{2/7}+\log\log\left(\frac{1}{\epsilon}\right)\big)$ of \eqref{bound:lip-uc}, which matches the lower bound given in \cite{arjevani2017oracle}. According to \eqref{eq:sublinear-res}, when $p=\nu=1, s=3$, with $q=2$, our algorithm has the iteration complexity $O\big(\sqrt{\frac{L}{\sigma}}\left(\frac{\sigma}{\epsilon}\right)^{\frac{1}{6}}\big)$, which may be of independent interest for solving the cubic regularized Newton step \cite{nesterov2006cubic} by gradient descent methods \cite{carmon2016gradient}. 

\subsection{Our Approach}\label{sec:approach}
In this paper, instead of directly designing an algorithm and then analyzing its iteration complexity, we consider a different paradigm to make our approach and algorithm more intuitive and explainable. 
The paradigm is inspired by the unified theory for first-order algorithms \cite{diakonikolas2019approximate} and the continuous-time interpretations of Nesterov's acceleration \cite{su2014differential,krichene2015accelerated,krichene2016adaptive,wibisono2016variational}. Our approach to the algorithmic design is based on an idealized but impractical algorithm called  \emph{vanilla proximal method (VPM)}, introduced in Section \ref{sec:vpm}. A continuous-time approximation to the VPM and a discrete-time approximation to the VPM will lead us to the final implementable algorithm with desired convergence rates. 
 
The VPM aims to solve a regularized program of the original one with an arbitrary convergence rate depending on parameters of our choice. However, the VPM serves more as an ideal target and is itself computationally infeasible to realize. We show that,in Section \ref{sec:cont}, to overcome the computational hurdle, one can instead solve a continuous-time \emph{convex approximation} to the VPM. An accelerated continuous-time  dynamics can be derived simply as sufficient conditions to ensure that solution to the approximate convex program achieves the same convergence rate as the original VPM. Such point of view unifies the existing continuous-time accelerated dynamics introduced in \cite{su2014differential}, \cite{krichene2015accelerated} and  \cite{wibisono2016variational} and severs as an arguably better guideline for the design of practical algorithms in the discrete setting.

In practice, to realize the desired accelerated dynamics, we need to know how to implement them in the discrete setting as an iterative algorithm. To this end, we need to consider a discrete-time {\em convex approximation} to the VPM. However, as we will see in Section \ref{sec:dis}, in order for the discrete-time approximation to achieve the same convergence rate as the continuous dynamics, we must solve a fixed-point problem which itself is computationally infeasible (if not impossible) in practice. To circumvent this difficulty, we propose to solve the fixed-point problem {\em approximately} by solving a  \emph{smooth approximation} to the VPM which becomes a tractable problem. Finally, by combing the convex approximation and the smooth approximation to the VPM, we propose an implementable discrete-time accelerated algorithm which achieves the optimal iteration complexity given in  \eqref{eq:rate} for the minimization of convex functions  with $(p, \nu, L)$-H\"{o}lder continuous derivatives (for $p\in\{1,2\ldots\}, \nu\in[0,1], L>0$).

\section{Preliminaries}\label{sec:prelim}
Before we proceed, we first introduce some notations. Let $:=$  denote a definition. Let $[n]$ denote the set $\{1,2,\ldots, n\}$. For $p=\{1,2,\ldots\}$, let $p! := 1\times2\times\cdots\times p$ with $0!:=1$. 
Let $\|\cdot\|$ denote a norm of vectors and $\|\cdot\|_*$ denote the dual norm of $\|\cdot\|$. For $x\in\bbR^d$ and $q\ge 1$, Let $\|x\|_q := (\sum_{i=1}^d |x_i|^q)^{\frac{1}{q}} $. { For a matrix $B\in\bbR^{d\times d}$ and $p,q\ge1$, denote the operator norm $\|B\|_{p,q}:=\max_{x\in\bbR^d}\{\|Bx\|_q: \|x\|_p\le 1\}$.}
By a little abuse of notation, for a convex function $f(x)$ defined on $\bbR^d$, let $\nabla f(x)$ denote the gradient at $x$ or one point in the subgradient set $\partial f(x)$. For a function $f(x;y)$,  $x$ denotes the variable of $f(x;y)$, $y$ denotes the parameter of $f(x;y)$ and $\nabla f(x;y)$ denotes the gradient or one point in the subgradient set $\partial f(x;y)$  $w.r.t.$ $x$. 

Similar to the notations in \cite{nesterov2018implementable}, for $p\in\{1,2,\ldots\}$, we use $\nabla^p f(x)[h_1, h_2, \ldots, h_p]$ to denote the directional derivative of a function $f$ at $x$ along the directions $h_i\in \bbR^d, i=1,2\ldots, p$. Then $\nabla^p f(x)[\cdot]$ is a symmetric $p$-linear form and its operator norm $w.r.t.$ a norm $\|\cdot\|$ is defined as
 \begin{eqnarray}
 \|\nabla^p f(x)\|_{*}:=\max_{y_1, y_2, \ldots,y_p}\left\{\nabla^p f(x)[y_1,\ldots, y_p]:\|y_i\|\le 1, i=1,2,\ldots, p\right\}.\label{eq:operator-norm}
 \end{eqnarray}

\begin{definition}[Strictly, Uniformly, or Strongly Convex]\label{def:convex}
We say a continuous function $f(x)$ is convex on $\bbR^d$, if $\;\forall x, y\in \bbR^d$, one has 
\begin{equation}
f(y)\ge f(x)+\langle \nabla f(x), y-x\rangle;	\label{eq:convex-def}
\end{equation}

$f(x)$ is strictly convex on $\bbR^d$, if the equality sign in \eqref{eq:convex-def} holds if and only if $x=y$;

$f(x)$ is $(s, \sigma)$-uniformly convex on $\bbR^d$ $w.r.t.$ a norm $\|\cdot\|$, if $\;\forall x,y\in \bbR^d$, one has 
\begin{eqnarray}
f(y)\ge f(x)+\langle \nabla f(x), y-x\rangle + \frac{\sigma}{s}\|y-x\|^s,\label{eq:uniform-def}	
\end{eqnarray}
where $s\ge 2$ is the order of uniform convexity and $\sigma\ge 0$ the constant of uniform convexity;

$f(x)$ is $\sigma$-strongly convex on $\bbR^d$ $w.r.t.$ $\|\cdot\|$, if $f(x)$ is $(2, \sigma)$-uniformly convex on $\bbR^d$ $w.r.t.$ $\|\cdot\|$.
\end{definition}
In Definition \ref{def:convex}, uniform convexity can be viewed as an extension of the better known concept of strong convexity. 
Example \ref{exm:1} gives two cases of  uniform convexity.
\begin{example}[Uniform Convexity]\label{exm:1}
$\frac{1}{2}\|x\|_q^2 (1<q\le 2)$ is $(2,q-1)$-uniformly convex on $\bbR^d$ $w.r.t.$  $\|\cdot\|_q$  \cite{ball1994sharp};  $\frac{1}{q}\|x\|_2^q\;(q\ge 2)$ is $(q, 2^{2-q})$-uniformly convex on $\bbR^d$ $w.r.t.$ $\|\cdot\|_2$ \cite{nesterov2008accelerating}.
\end{example}

Starting from the work of \cite{nesterov2015universal}, an increasing interest is to replace the Lipschitz continuity assumption by the H\"{o}lder continuity assumption \cite{yurtsever2015universal,roulet2017sharpness,nesterov2018complexity,grapiglia2019tensor} and to propose universal algorithms in the sense that the convergence of algorithms can optimally adapt to the H\"{o}lder parameter. 
\cite{nesterov2015universal,yurtsever2015universal} have considered first-order algorithms with H\"{o}lder continuous gradients $w.r.t.$ $\|\cdot\|_2$; \cite{grapiglia2019accelerated} has proposed cubic regularized Newton methods for minimizing functions with H\"{o}lder continuous Hessians $w.r.t.$ $\|\cdot\|_2$; \cite{cartis2019universal,grapiglia2019tensor} have considered tensor methods for minimizing functions with $p$-th  H\"{o}lder continuous derivatives $w.r.t.$ $\|\cdot\|_2$ ($p\in\{2,3,\ldots\}$). In this paper, we extend the definition of H\"{o}lder continuous derivatives $w.r.t.$ any norm $\|\cdot\|$, including non-Euclidean norms. Our analysis and results will be applicable to this general setting.

\begin{definition}[H\"{o}lder Continuous Derivative]\label{def:holder-deri}
We say a function $f(x)$ on $\bbR^d$ has $(p, \nu, L)$-H\"{o}lder continuous derivatives $w.r.t.$ $\|\cdot\|$, if $\;\forall x, y\in \bbR^d$,  one has

\begin{eqnarray}
\frac{1}{(p-1)!}\|\nabla^p f(x) - \nabla^p f(y)\|_{*} \le  L\|x-y\|^{\nu}, 
\end{eqnarray}
where $p\in\{1, 2,3,\ldots\}$ denotes the order of derivative, $0\le\nu\le 1$ denotes the H\"{o}lder parameter and $L>0$ is the constant of smoothness.

$f(x)$ is said to have $(p, L)$-Lipschitz continuous derivatives on $\bbR^d$ $w.r.t.$ $\|\cdot\|$  if $f(x)$ has $(p, 1, L)$-H\"{o}lder continuous derivatives on $\bbR^d$ $w.r.t.$ $\|\cdot\|$.
\end{definition}
In Definition \ref{def:holder-deri}, we unify the definition of H\"{o}lder continuous gradients $(i.e., p =1)$ and high-order H\"{o}lder continuous derivatives $(i.e., p \in\{2,3,\ldots\})$. For $p=1$, $\|\cdot\|_*$ denotes the dual norm of $\|\cdot\|$; for $p\in\{2,3,\ldots\}$, $\|\cdot\|_*$ denotes the operator norm of tensor of $p$-th order $w.r.t.$ $\|\cdot\|$, which is defined by  \eqref{eq:operator-norm}.

{
\begin{example}[Non-Euclidean High-order Smoothness]\label{exam:smooth}
Consider the objective function $f(x):= \frac{1}{n} \sum_{j=1}^n \log(1+\exp(-\bar{b}_j \bar{a}_j^T x))$ for logistic regression, where  $j\in[n], \bar{a}_j\in\bbR^d, \bar{b}_j\in\{1,-1\}.$ Denote $B:=\frac{1}{n}\sum_{j=1}^n \bar{a}_j \bar{a}_j^T.$ For $1\le p\le 2$ and $q$ satisfying $1/p + 1/q = 1$, let $\|\nabla^s f(x)\|_{q}$ denote the operator norm of $\nabla^s f(x) (s=2,3)$  in \eqref{eq:operator-norm} $w.r.t.$ the vector norm $\|\cdot\|_p$. Then we have
\begin{equation}
\|\nabla^2 f(x) -   \nabla^2 f(y)\|_{q} \le \|B\|_{p,q}\max_{j\in[n]}\|\bar{a}_j\|_q^{\nu}\cdot \|x - y\|^{\nu}_p.  
\end{equation}
\end{example}

\begin{proof}
See Section \ref{sec:exam:smooth}.
\end{proof}

}

In the paper, we mainly consider the problem of optimizing a composite convex function of the form
\begin{eqnarray}
\min_{x\in \bbR^d} f(x):=g(x)+l(x),\label{eq:prob}
\end{eqnarray}
where $g(x)$ is a closed proper convex function and $l(x)$ is a simple convex but maybe non-smooth function. We consider the case when $g(x)$ has $(p, \nu, L)$-H\"{o}lder continuous derivatives, for all $x, y\in\bbR^d$. Then we can define the following two auxiliary functions that approximate $f(x)$:
\begin{eqnarray}
\hat{f}(x;y)&:=&g(y)+\langle \nabla g(y), x-y\rangle + l(x), \label{eq:hat-f}\\
\tilde{f}(x;y)&:=& g(y)+ \sum_{i=1}^p\frac{1}{i!}\nabla^i g(y)[x-y]^i + l(x). \label{eq:tilde-f}
\end{eqnarray}
Then $\hat{f}(x;y)$ is a lower-bound convex approximation to $f(x)$ for any parameter $y\in\bbR^d$. $\tilde{f}(x;y)$ gives a high-order smooth approximation to $f(x)$ for any parameter $y\in\bbR^d$. Or more formally, we have:
\begin{lemma}\label{lem:hat-tilde-f}
 If $g(x)$ and $ l(x)$ are convex, and $g(x)$ has $(p, \nu, L)$-H\"{o}lder continuous derivatives, for all $x, y\in\bbR^d$, then we have
\begin{eqnarray}
\hat{f}(x;y)&\le& f(x)\label{eq:hat-prop},\\
|f(x) - \tilde{f}(x;y)|&\le&\frac{L}{p}\|x-y\|^{p+\nu},\label{eq:tilde-prop-1}\\
\|\nabla f(x) -\nabla \tilde{f}(x;y)\|_*&\le&{L}\|x-y\|^{p+\nu-1}.\label{eq:tilde-prop-2}
\end{eqnarray}
\end{lemma}
\begin{proof}
See Section \ref{sub:lem:hat-tilde-f}. 
\end{proof}
In \eqref{eq:hat-f} and \eqref{eq:tilde-f}, we do not linearize the term $l(x)$ which may be nonsmooth. 
Because of \eqref{eq:hat-prop}, in this paper, $\hat{f}(x;y)$ is viewed as a lower-bound convex approximation to $f(x)$ for any parameter $y\in\bbR^d$. $\tilde{f}(x;y)$ satisfies \eqref{eq:tilde-prop-1} and \eqref{eq:tilde-prop-2}, and gives a high-order smooth approximation to $f(x)$ for any parameter $y\in\bbR^d$.
In our analysis, the convexity and smoothness assumptions are used only by the two inequalities \eqref{eq:hat-prop} and \eqref{eq:tilde-prop-2}, which allow a unified treatment for the smooth  and the composite convex settings (with or without the term $l(x)$). Meanwhile, because we only need the property  \eqref{eq:tilde-prop-2} of high-order smoothness, it implies that in convex optimization, the high-order smoothness is mainly used to give a more accurate estimation of the ``implicit gradient'' $\nabla f(x).$

Finally, we give two inequalities in Lemma \ref{lem:seq} which will be used in our analysis. 
\begin{lemma}\label{lem:seq}
For a  sequence $\{b_k\}_{k\ge 0}$ with $b_0=0$ and $b_k>0\; (k\ge 1)$. Then for $\rho \ge 1$ and $C>0$, if $\;\forall k\ge1$,
\begin{eqnarray}
(b_k-b_{k-1})^{\rho}\ge C b_k^{\rho-1}, \label{eq:seq-1}
\end{eqnarray}
then we have
\begin{eqnarray}
b_k \ge C\left(\frac{k}{\rho}\right)^{\rho}. \label{eq:seq-2}
\end{eqnarray}
Meanwhile, for $\rho \ge 1, \delta >0$ and $C>0$, 
if $\;\forall k \ge 1,$
\begin{eqnarray}
\sum_{i=1}^k \left(\frac{b_i^{\rho-1}}{(b_i-b_{i-1})^{\rho}}\right)^\delta \le C, \label{eq:seq-3}
\end{eqnarray}
then we have
\begin{eqnarray}
b_k\ge C^{-\frac{1}{\delta}}\left(\frac{k}{\rho}\right)^{\rho+\frac{1}{\delta}}.  \label{eq:seq-4}
\end{eqnarray}
\end{lemma}

\begin{proof}
See Section \ref{sub:lem:seq}. 
\end{proof}

\section{A Vanilla Proximal Method}\label{sec:vpm}
Let us start our study by considering a composite convex optimization problem in \eqref{eq:prob}.
In the following discussion, we assume that $x^*$ is a minimizer of $f(x)$ on $\bbR^d$. To design an acceleration algorithm to minimize $f(x)$, we first introduce a so-called \emph{vanilla proximal method (VPM)}, that considers to minimize an auxiliary function $\psi^{\text{vpm}}(x):=A f(x) + h(x; x_0)$ as in Algorithm \ref{alg:vpm}.
\begin{algorithm} [!ht]
\caption{Vanilla Proximal Method (VPM)}\label{alg:vpm}
\begin{algorithmic}[1]
\STATE \textbf{Input:} an initialized point $x_0\in \bbR^d$, a positive scalar  $A>0$. 
\STATE Find a $z\in\bbR^d$ that satisfies 
\begin{equation}
z := \argmin_{x\in \bbR^d}\big\{\psi^{\text{vpm}}(x):=A f(x) + h(x; x_0) \big\}. \label{eq:vpm}
\end{equation}
\RETURN $z$.
\end{algorithmic}	
\end{algorithm}

In the auxiliary objective $\psi^{\text{vpm}}(x)$, the proxy term $h(x;x_0)$ typically should satisfy the following assumption: 
\begin{assumption}\label{ass:h-vpm}
For all $x,x_0\in \bbR^d$, $h(x;x_0)\ge 0$ with $h(x;x_0)=0$ if and only if $x=x_0$. Meanwhile, $h(x;x_0)$ is strictly convex on $\bbR^d$. 
\end{assumption}
Therefore, in the VPM, \eqref{eq:vpm} is a strictly convex program and thus there exists a unique minimizer $z$. By using only the optimality condition of \eqref{eq:vpm}, we can characterize the ``convergence rate'' of the VPM as below.

\begin{theorem}\label{thm:vpm}
The solution $z$ generated by Algorithm \ref{alg:vpm} satisfies
\begin{eqnarray}
f(z) - f(x^*)\le \frac{h(x^*; x_0)}{A}. 
\end{eqnarray}
\end{theorem}
\begin{proof}%
By the definition of $\psi^{\text{vpm}}(x)$ in \eqref{eq:vpm}, one has 
\begin{align}
\min_{x\in \bbR^d} {\psi^{\text{vpm}}(x)}\le A f(x^*) + h(x^*; x_0).\label{eq:vpm-upper}
\end{align}
Then by the optimality condition of $z$ and the nonnegativity of $h(x; x_0)$, one has
\begin{align}
A f(z){\le} A f(z) + h(z; x_0)=\min_{x\in \bbR^d} {\psi^{\text{vpm}}(x)}. \label{eq:vpm-lower}
\end{align}
By the upper bound of $\min_{x\in \bbR^d} {\psi^{\text{vpm}}(x)}$ in \eqref{eq:vpm-upper} and lower bound of $\min_{x\in \bbR^d} {\psi^{\text{vpm}}(x)}$ in \eqref{eq:vpm-lower}, after a simple rearrangement, Theorem \ref{thm:vpm} is proved.

\end{proof}

According to Theorem \ref{thm:vpm}, the VPM may converge with any convergence rate if $A$ is chosen to a large enough value. In fact we do not need any extra assumption on $f(x)$ in the proof of Theorem \ref{thm:vpm}, except that the optimal solution  $z$ exists. Although solving the subproblem \eqref{eq:vpm} is impractical in general, it provides us a good starting point to design practical algorithms: by making certain assumptions on the objective function $f(x)$ and the proxy function $h(x;x_0)$, it is possible to achieve or approach the convergence rate of the VPM by solving a tractable approximation to \eqref{eq:vpm}. {The ideal subproblem \eqref{eq:vpm} is non-iterative and does not define an optimization path. Nevertheless, when we consider a tractable approximation to \eqref{eq:vpm}, the approximation can be defined by a continuous-time or discrete-time optimization path with dependence on previous states or iterations.    
}
As we will see, a continuous approximation results in a continuous-time accelerated  dynamic system in Section \ref{sec:cont}, while a discrete-time approximation results in a discrete-time accelerated algorithm in Section \ref{sec:dis}.

\begin{remark}
The proposed VPM is similar to the proximal point algorithm (PPA) \cite{parikh2014proximal} which performs \begin{equation}
x_{k+1} := \argmin_{x\in \bbR^d}\Big\{ a_k f(x)+\frac{1}{2}\|x-x_k\|_2^2\Big\}\label{eq:ppa}	
\end{equation}
along iterations, where $a_k>0$.   The difference between VPM and PPA is that VPM is not an iterative algorithm and has a convergence rate only depending on the parameter $A$ we choose. If we set $h(x;x_0) = \frac{1}{2}\|x-x_0\|_2^2$, the per-iteration costs of VPM and PPA are comparable.
\end{remark}
\medskip

\begin{remark}
The subproblem \eqref{eq:ppa} of PPA is often computationally infeasible in practice.
By considering tractable inexact versions of PPA with the concept  ``$\epsilon$-subdifferential'', \cite{monteiro2013accelerated} has proposed the unified A-HPE framework for convex optimization. One difference between our framework and A-HPE is that ours extends from the non-iterative VPM and therefore can unify both continuous-time accelerated dynamics and discrete-time accelerated algorithms. Meanwhile, we consider a general proxy function $h(x;x_0)$ rather than the Euclidean norm square $\frac{1}{2}\|x-x_0\|_2^2$, thus our framework can be generalized to the non-Euclidean setting such as $\frac{1}{2}\|x-x_0\|_q^2 \;(1<q\le 2)$ and the $q$-th power of Euclidean norm $\frac{1}{q}\|x-x_0\|_2^q\;(q\ge 2)$. 
\end{remark}

{\rc

}

\section{Continuous-time Accelerated Descent Dynamics}\label{sec:cont}
The subproblem \eqref{eq:vpm} in the VPM is merely conceptual as it is almost as difficult as minimizing the original function. Nevertheless, if $f(x)$ is convex, one can always seek  more tractable approximations.
From an acceleration perspective, the convex approximation $\hat{f}(x; y)$ in Lemma \ref{lem:hat-tilde-f} gives a lower bound for $f(x)$ at the  state $y$. The minimizer of $\hat{f}(x; y)$ would suggest an aggressive direction and step for the next iterate to go to. However, for such iterates not to diverge too far from the landscape of $f(x)$, we also need a good upper bound. A basic idea is that up to time $t$, we have already traversed a path $x_{\tau}, \tau \in [0, t)$ over the landscape of $f(x)$. We could potentially use all the lower-bounds $\hat{f}(x; x_{\tau})$ of $f(x)$ to construct a good upper bound to guide the next step. The simplest possible form for such an upper bound we could consider is a superposition (or an integral) of these lower bounds:
$$\psi_t^{\text{cont}}(x):=\int_{0}^t a_{\tau} \hat{f}(x; x_{\tau}) \d\tau+h(x; x_0),$$
where $a_\tau$ is a properly chosen weight function of $\tau$ and $h(x; x_0)$ is a strictly convex term to bound the function $\psi_t^{\text{cont}}$ from below in case $\hat{f}(x; x_{\tau})$ are not.

Therefore, to guide the descent trajectory, we can consider solving an approximate problem of \eqref{eq:vpm} as follows
\begin{eqnarray}
z_t := \argmin_{x\in \bbR^d}\left\{\psi_t^{\text{cont}}(x):=\int_{0}^t a_{\tau} \hat{f}(x; x_{\tau}) \d\tau+h(x; x_0)\right\}, \label{eq:cont}
\end{eqnarray}
where  $\forall\; 0< \tau\le t, a_{\tau}> 0$ and satisfies $\int_{0}^t a_{\tau}\d\tau = A_t$ with $a_0 = A_0=0$  and  $\{x_{\tau}\}_{0\le \tau\le t}$ is the optimization path and its relationship with $z_t$ will be determined soon.

In this section, our main goal is to show that the widely studied continuous-time accelerated dynamics arise from a sufficient condition that allows \eqref{eq:cont} to achieve the same convergence rate as the original VPM. First, a upper bound of $\min_{x\in \bbR^d}\psi_t^{\text{cont}}(x)$ is given in Lemma \ref{lem:cont-upper}.
\begin{lemma}\label{lem:cont-upper}
For all $t \ge 0$, we have $\min_{x\in \bbR^d}\psi_t^{\rm{cont}}(x)\le A_t f(x^*) + h(x^*; x_0)$. 
\end{lemma} 
\begin{proof}
See Section \ref{sub:lem:cont-upper}
\end{proof}
Lemma \ref{lem:cont-upper} is an extension of the upper bound \eqref{eq:vpm-upper} of $\psi^{\rm{vpm}}$, which follows trivially from  Lemma \ref{lem:hat-tilde-f}. In other words, Lemma \ref{lem:cont-upper} provides a lower bound of $A_t f(x^*)$.

\begin{lemma}\label{lem:cont}
 $\forall t\ge 0$, we have $A_t f(x_t) \le \min_{x\in \bbR^d} \psi_t^{\rm{cont}}(x) + \int_{0}^t \langle \nabla f(x_{\tau}), A_{\tau}\dot{x}_{\tau}- a_{\tau}(z_{\tau} - x_{\tau})\rangle \d\tau$.
\end{lemma}
\begin{proof}
See Section \ref{sub:lem:cont}
\end{proof}
Essentially, Lemma \ref{lem:cont} says that the lower bound \eqref{eq:vpm-lower} of $\psi^{\rm{vpm}}$ can be extended to $\psi_t^{\rm{cont}}$, at least approximately. We would like to make this approximation as close as possible and establish $\min_{x\in\bbR^d}\psi_t^{\text{cont}}(x)$ as an upper bound of $A_t f(x_t)$, at least along certain path by our choice. To this end, based on Lemmas \ref{lem:cont-upper} and \ref{lem:cont}, we have the following theorem.
\begin{theorem}[Continuous-Time VPM]\label{thm:cont}
If the continuous-time trajectories $\{x_t\}_{t\ge 0}$ and $\{z_t\}_{t\ge 0}$ are evolved according to the dynamics:
\begin{equation}
\begin{cases}
A_{t}\dot{x}_{t} = a_{t}(z_{t}-x_{t}), \\
z_t = \argmin_{x\in \bbR^d}\left\{\int_{0}^t a_{\tau} \hat{f}(x; x_{\tau}) \d\tau+h(x; x_0)\right\},
\end{cases}\label{eq:dynamic}
\end{equation}
where  $\forall \;0<\tau\le t$, $a_{\tau}> 0$, $\int_{0}^t a_{\tau} \d\tau = A_t $, and $a_0=A_0=0$,
then for all $t>0$, one has
\begin{align}
f(x_t) - f(x^*)\le \frac{h(x^*; x_0)}{A_t}.\label{eq:rate-con}
\end{align}
\end{theorem}
\begin{proof}
If $A_{t}\dot{x}_{t} = a_{t}(z_{t}-x_{t})$, then from Lemma \ref{lem:cont}, one has 
$$A_t f(x_t) \le \min_{x\in \bbR^d} \psi_t^{\rm{cont}}(x).$$
Combining Lemma \ref{lem:cont-upper}, we have \eqref{eq:rate-con}.
\end{proof}
\medskip

In  Theorem \ref{thm:cont}, \eqref{eq:dynamic} does not specify any concrete values or forms for $a_t$ and $A_t$, except the condition $a_{\tau}> 0$, $\int_{0}^t a_{\tau} \d\tau = A_t  \;(0<\tau\le t)$ and $a_0=A_0=0$;\footnote{Theoretically $A_t$ should be chosen such that the differential equation to have a unique solution.} meanwhile it does not specify any concrete form for $h(x;x_0)$. As result, by instantiating the dynamical system \eqref{eq:dynamic}, one may obtain all the ODEs previously introduced and studied in the literature \cite{su2014differential,krichene2015accelerated,krichene2016adaptive,wibisono2016variational}, respectively. We show a few examples below:

\begin{example}\label{exm:2}
If  the component of \eqref{eq:prob} $l(x):=0$ and $ h(x;x_0) := \frac{1}{2}\|x-x_0\|_2^2$, then \eqref{eq:dynamic} is equivalent to
\begin{equation}
\ddot{x}_t+\frac{a_t}{A_t}\left(\frac{\d}{\d t}\left(\frac{A_t}{a_t}\right)+1\right)\dot{x}_t + \frac{a_t^2}{ A_t}\nabla f(x_t) = 0, \label{eq:Eucli}
\end{equation}
where  $\forall\; 0<\tau\le t$, $ a_{\tau}> 0$, $A_t=\int_{0}^t a_{\tau} \d\tau, a_0=A_0=0$. By setting $A_t := \frac{1}{4}t^2$ and $A_t := \frac{1}{p^2}t^p\;(p>2)$ respectively, then we recover the ODE in \cite{su2014differential} and the ODE  under the Euclidean norm setting  in \cite{wibisono2016variational}.
\end{example}

\begin{remark}
In Example \ref{exm:2}, if $l(x)$ is the indicator function of a closed convex set and $h(x;x_0)$ is chosen as the Bregman divergence of a strictly convex function, then we may recover the formulation of accelerated mirror descent dynamic \cite{krichene2015accelerated} and the Euler-Lagrange equation \cite{wibisono2016variational}. 
\end{remark}

\begin{remark}
Although we have derived the dynamics \eqref{eq:dynamic} from a different perspective, it should be noted that the dynamical system \eqref{eq:dynamic} is an extension and refinement to the ODE derived by the ``approximate duality gap technique (ADGT)'' \cite{diakonikolas2019approximate}. The main difference is that instead of giving a upper bound of $f(x_t)$ and a lower bound of $f(x^*)$, we give a upper bound of $A_tf(x_t)$ and a lower bound of $A_tf(x^*)$. This modification allows us to set $A_0=0$ rather than $A_0>0$, and thus the initialization expression about $A_0$ can be removed. 
Such a modification simplifies future derivation and analysis greatly.
\end{remark}

Compared to the VPM, the continuous-time accelerated dynamics must satisfy an extra ODE, which can be viewed as an additional cost associated with the  continuous-time approximation. As we see in Theorem \ref{thm:cont}, the optimization path $\{x_t\}_{t\ge 0}$ can have the same property of the VPM and one can obtain an arbitrarily fast convergence rate if $A_t$ is chosen to be large enough. However, if a discrete-time approximation is used to implement and approximate the VPM, it is in general difficult to retain the same rate, which will be discussed carefully in the next Section \ref{sec:dis}.

\section{Discrete-time Accelerated Descent Algorithm}\label{sec:dis}

In order to achieve the same convergence rate of the VPM, the continuous-time approximation needs the extra ODE condition in \eqref{eq:dynamic}, which is reasonable to assume in the continuous setting. In the discrete-time setting, if all other conditions remain unchanged, except that we replace the weighted continuous-time approximation \eqref{eq:cont} by a weighted discrete-time counterpart, one may see that the ODE will be replaced by a condition that requires us to find a solution to a \emph{fixed-point} problem (which will be clear in Lemma \ref{lem:dis-lower}). Unfortunately, directly solving this fixed-point problem is computationally infeasible in practice. To remedy this difficulty, we employ  a stronger assumption for the proxy function $h(x;x_0)$ such that it can introduce an extra term $\frac{\gamma}{q}\|y-x\|^q$ as follows.
\begin{assumption}\label{ass:h-dis}
For all $x,x_0\in \bbR^d$, $h(x;x_0)\ge 0$ with $h(x;x_0)=0$ if and only if $x=x_0$. Meanwhile $h(x;x_0)$ is $(q, \gamma)$-uniformly convex  $w.r.t.$ a  norm $\|\cdot\|$ where $q\ge 2, \gamma> 0$, $i.e.,$ $\forall x,y\in\bbR^d$, one has
\begin{align}
h(y;x_0)\ge h(x;x_0) + \langle \nabla h(x;x_0), y-x\rangle + \frac{\gamma}{q}\|y-x\|^q.
\end{align}
\end{assumption}

\begin{example}
By setting $h(x;x_0):=\frac{1}{2}\|x-x_0\|_p^2\;(1<p\le 2)$ or the Bregman divergence of $\frac{1}{2}\|x-x_0\|_p^2$, then $h(x;x_0)$ satisfies Assumption \ref{ass:h-dis} $w.r.t.$ $\|\cdot\|_p$ with order $q=2$ and constant $\gamma = p-1$ \cite{ball1994sharp}; by setting $h(x;x_0):=\frac{1}{p}\|x-x_0\|_2^p \;(p\ge 2)$, then $h(x;x_0)$ satisfies Assumption \ref{ass:h-dis} $w.r.t.$ $\|\cdot\|_2$ with order $q=p$ and constant $\gamma = 2^{2-p}$ \cite{nesterov2008accelerating}.
\end{example}
Meanwhile, we also need that $\frac{1}{q}\|\cdot\|^q$ is $(q, \beta)$-uniformly convex.          
\begin{assumption}\label{ass:norm}
For the norm $\| \cdot\|$, we assume that $\frac{1}{q}\|x\|^q$ is $(q, \beta)$-uniformly convex $w.r.t.$ $\|\cdot\|$, where $q\ge 2, \beta>0$,  $i.e.,$ 
\begin{align}
\frac{1}{q}\|y\|^q\ge\frac{1}{q} \|x\|^q + \left\langle \nabla\frac{1}{q} \|x\|^q, y-x\right\rangle + \frac{\beta}{q}\|y-x\|^q.
\end{align}
\end{assumption}

In order to find a good approximate solution of our problem in a computationally efficient way, the smooth component $g(x)$ of $f(x)$ in \eqref{eq:prob} should satisfy the following. 
\begin{assumption}\label{ass:f-smooth}
$g(x)$ has $(p, \nu, L)$-H\"{o}lder continuous derivatives, where $p\in \{1, 2, \cdots\}$, $\nu\in[0,1], L>0$.
\end{assumption}
In Assumptions \ref{ass:h-dis} to \ref{ass:f-smooth}, for practical concerns and  technical reasons, in the following discussion, we will assume that $p+\nu \ge 2$ and $q\in[2, p+\nu]$. ($p+\nu\ge 2$ means that in our setting if $p=1$, then $\nu=1$.)

Based on Assumptions \ref{ass:h-dis}-\ref{ass:f-smooth},  
in this section, similar to the continuous-time approximation in Section \ref{sec:cont}, we consider a weighted discrete-time convex approximation of \eqref{eq:vpm}: for $k\ge 0$, 
\begin{eqnarray}
z_k := \argmin_{x\in \bbR^d}\left\{\psi_k^{\text{dis}}(x):=\sum_{i=1}^{k} a_{i} \hat{f}(x; x_{i}) + h(x;x_0)\right\}, \label{eq:dis}
\end{eqnarray}
where we assume that $\forall 1\le i\le k, a_{i}> 0$,   $A_i:=\sum_{j=1}^i a_j$ and  $a_0=A_0=0$, $h(x;x_0)$ satisfies Assumption \ref{ass:h-dis}, and $\hat{f}(x; x_{i})$ is defined in Lemma \ref{lem:hat-tilde-f}. Meanwhile, in \eqref{eq:dis}, when $k=0$, we let $\psi_0^{\text{dis}}(x) =  h(x;x_0)$ and thus $z_0=\argmin_{x\in\bbR^d} h(x;x_0) = x_0$.
Then we motivate the discrete-time algorithm by analyzing the conditions needed to emulate the same rate of the VPM.
First, a upper bound of $\min_{x\in \bbR^d}\psi_k^{\text{dis}}(x)$ is given in Lemma \ref{lem:dis-upper}.
\begin{lemma}\label{lem:dis-upper}
For all $k\ge 0$, one has $\min_{x\in \bbR^d}\psi_k^{{\rm{dis}}}(x)\le A_k f(x^*) + h(x^*;x_0)$. 
\end{lemma}
\begin{proof}
see Section \ref{subsec:lem:dis-upper}.
\end{proof}

Then in Lemma \ref{lem:dis-lower} below, we show how the lower bound \eqref{eq:vpm-lower} of $\psi^{\rm{vpm}}$ can be extended to the discrete case $\psi_k^{\rm{dis}}$ with some extra terms.

\begin{lemma}\label{lem:dis-lower}
For $i\ge 1$, let $E_{i}:=A_{i}  \left\langle \nabla f(x_{i}), x_{i}-\frac{a_{i}}{A_{i}}z_{i}-\frac{A_{i-1}}{A_{i}}x_{i-1}\right\rangle - \frac{\gamma}{q}\|z_{i}-z_{i-1}\|^q$. 
Then for all $k\ge 1$, one has 
\begin{align}
&A_k f(x_k)-  \psi_k^{\rm{dis}}(z_k)\le \sum_{i=1}^{k} E_{i}.\label{eq:dis-lower}
\end{align}
\end{lemma}
\begin{proof}
	See Section \ref{subsec:lem:dis-lower}.
\end{proof}
In Lemma \ref{lem:dis-lower}, the extra negative term $- \frac{\gamma }{q}\|z_{i}-z_{i-1}\|^q$ in $E_i$ is from the uniform convexity of $h(x;x_0)$. If $h(x;x_0)$ is only convex (i.e. $\gamma =0$), this negative term does not exist and thus a sufficient condition for $E_{i}\le 0$ is: 
\begin{eqnarray}
 x_{i}=\frac{a_{i}}{A_{i}}z_{i}+\frac{A_{i-1}}{A_{i}}x_{i-1}, \quad \forall 1\le i\le k. \label{eq:x-i-1}
\end{eqnarray}
By \eqref{eq:dis}, $z_{i}$ is a function of $x_{i}$. Therefore finding $x_{i}$ to satisfy \eqref{eq:x-i-1} is reduced to a fixed-point problem (so is it for $z_{i}$). It is computationally infeasible (if not impossible) to find an exact solution to this problem in general. Nevertheless, if $h(x;x_0)$ satisfies Assumption \ref{ass:h-dis}, the term
$E_{i}$ contains a negative term $- \frac{\gamma}{q}\|z_{i}-z_{i-1}\|^q$. So there is hope that an approximate solution to the fixed-point problem \eqref{eq:x-i-1} can still make $E_{i}\le 0$.

To approximately solve the fixed-point problem, for convenient analysis, inspired by \cite{hairer1987solving,diakonikolas2017accelerated}, we define a pair $(\hat{x}_{i-1}, \hat{z}_{i})$ such that 
\begin{equation}
\begin{cases}
\hat{x}_{i-1} := \frac{a_{i}}{A_{i}}z_{i-1}+\frac{A_{i-1}}{A_{i}}x_{i-1}, \\
\hat{z}_{i}:=\frac{A_{i}}{a_{i}}x_{i} - \frac{A_{i-1}}{a_{i}}x_{i-1}.
\end{cases}\label{eq:two-fixed}
\end{equation}
By \eqref{eq:two-fixed}, we have $x_{i} = \frac{a_{i}}{A_{i}}\hat{z}_{i}+\frac{A_{i-1}}{A_{i}}x_{i-1}$. Therefore \eqref{eq:two-fixed} can be viewed as two-step fixed-point iterations for $x_{i}$ based on $\hat{x}_{i-1}$ and $\hat{z}_{i}$. Here $\hat{z}_{i}$ can be viewed as the best estimate of the desired fixed point $z_{i}$ based on the calculated $x_{i}$ in our algorithm. It is defined for convenience and will only be used in our analysis but not in the algorithm.

Based on the definition of $(\hat{x}_{i-1}, \hat{z}_{i})$, Assumption \ref{ass:norm}, 
and the definition of $E_{i}$ in Lemma \ref{lem:dis-lower},
we have Lemma \ref{lem:E-2}.

\begin{lemma}\label{lem:E-2}
For $i\ge 1$ and any $ \gamma^{\prime}_{i}\in (0, \gamma]$, we have
\begin{eqnarray}
E_{i} &\le& a_{i}\left\langle \nabla f(x_{i}) + \frac{\gamma^{\prime}_{i} A_{i}^{q-1}}{a_{i}^{q}}\nabla\frac{1}{q} \|x_{i} -\hat{x}_{i-1}\|^q,  \hat{z}_{i} - z_{i}\right\rangle \nonumber\\
&&-\gamma^{\prime}_{i}\left(
\frac{A_{i}^q}{qa_{i}^q}\|x_{i} -\hat{x}_{i-1}\|^q + \frac{\beta}{q}\|\hat{z}_{i} - z_{i}\|^q
  \right). \label{eq:E-2}
\end{eqnarray}
\end{lemma}

\begin{proof}
See Section \ref{subsec:lem:E-2}. 
\end{proof}
In Lemma \ref{lem:E-2}, we purposely introduce a new parameter $\gamma^{\prime}_{i}$, which as we will soon show, helps unify the four high-order instances \cite{monteiro2013accelerated, jiang2018optimal,bubeck2018near,jiang2018optimal} of the A-HPE framework. Meanwhile, because of the uniform convexity of $\frac{1}{q}\|\cdot\|^q$, the negative term $- \frac{\gamma}{q}\|z_{i}-z_{i-1}\|^q$ is reduced to two negative terms and a term of inner product.  

By Lemma \ref{lem:E-2}, if we can find $x_{i}$ such that  
\begin{eqnarray}
\nabla f(x_{i}) + \frac{\gamma^{\prime}_{i} A_{i}^{q-1}}{a_{i}^{q}}\nabla\frac{1}{q} \|x_{i} -\hat{x}_{i-1}\|^q=0,	\label{eq:E-3}
\end{eqnarray}
then we can ensure $E_{i}\le 0$. However the problem of finding $x_{i}$ satisfying \eqref{eq:E-3} is equivalent to solving the VPM problem exactly in \eqref{eq:vpm} with the settings $x_0 := \hat{x}_{i-1}, h(x;x_0) :=\frac{1}{q}\|x - \hat{x}_{i-1}\|^q,   A:=\frac{a_{i}^{q}}{\gamma^{\prime}_{i} A_{i}^{q-1}}$, which is computationally infeasible in general.  Fortunately, the two negative terms in \eqref{eq:E-2} may dominate small errors even if we solve the VPM problem \eqref{eq:E-3} inexactly. Hence we approximate the intermediate VPM problem \eqref{eq:E-3} by a smooth approximation $\tilde{f}(x_{i};\hat{x}_{i-1})$ using the fact
\begin{eqnarray}
\|\nabla f(x_{i}) - \nabla\tilde{f}(x_{i};\hat{x}_{i-1})\|_*\le L\|x_{i}-\hat{x}_{i-1}\|^{p+\nu-1}, \label{eq:E-4}
\end{eqnarray}
from Lemma \ref{lem:hat-tilde-f}.
Then by Lemmas \ref{lem:hat-tilde-f} and \ref{lem:E-2},  we have Lemma \ref{lem:E-3}.
\begin{lemma}\label{lem:E-3}
Denote $c_q := \left(\beta (q-1)^{1-q}\right)^{\frac{1}{q}}$ and $\lambda_{i}^{\prime}:= \frac{a_{i}^q}{c_q \gamma^{\prime}_{i} A_{i}^{q-1}}$. 
For $i\ge 1$, one has
\begin{eqnarray}
E_{i}&\le&\left(\left(L\lambda_{i}^{\prime}\|x_{i} -\hat{x}_{i-1}\|^{p+\nu-q} \right)^{\frac{q}{q-1}} -1\right)\frac{\gamma^{\prime}_{i} A_{i}^q}{q a_{i}^q}\|x_{i} -\hat{x}_{i-1}\|^q \nonumber\\
&&+a_{i}\left \langle \nabla\tilde{f}(x_{i};\hat{x}_{i-1})
+ \frac{\gamma^{\prime}_{i} A_{i}^{q-1}}{a_{i}^q}\nabla\frac{1}{q} \|x_{i} -\hat{x}_{i-1}\|^q ,  \hat{z}_{i} - z_{i}\right\rangle\label{eq:E-5}. 
\end{eqnarray}
\end{lemma}
\begin{proof}
Sec Section \ref{sec:E-3}.
\end{proof}

In Lemma \ref{lem:E-3}, by using the smooth approximation $\tilde{f}(x_{i};\hat{x}_{i-1})$, to ensure $E_{i}\le 0$,  the VPM problem \eqref{eq:E-3} is reduced to an easier smooth approximation problem
\begin{eqnarray}
\nabla\tilde{f}(x_{i};\hat{x}_{i-1})+ \frac{\gamma^{\prime}_{i} A_{i}^{q-1}}{a_{i}^q}\nabla\frac{1}{q} \|x_{i} -\hat{x}_{i-1}\|^q = 0
\label{eq:E-6-2}
\end{eqnarray}
with a cost that the introduced smooth approximation error should not go beyond the capability of the two negative terms \eqref{eq:E-2} to balancing our errors. As a result, in Lemma \ref{lem:E-3}, to ensure $E_{i}\le 0$, we also need the condition 
\begin{eqnarray}
L\lambda_{i}^{\prime}\|x_{i} -\hat{x}_{i-1}\|^{p+\nu-q}\le\theta_2\le 1 \label{eq:E-6-1}
\end{eqnarray}
is true, where $\theta_2\in(0,1]$ is a constant.

We here discuss the role of the parameter $\gamma^{\prime}_{i}$. So far our derivation works for any $\gamma^{\prime}_{i} \in(0,\gamma]$.
 A simple choice of $\gamma^{\prime}_{i}$ would be $\gamma^{\prime}_{i}:=\gamma$. Nevertheless,  
under the condition \eqref{eq:E-6-1}, for any $\alpha\in [0, 1]$, we could choose $\gamma^{\prime}_{i}$ to satisfy:
\begin{eqnarray}
\gamma^{\prime}_{i}:=\left( \frac{L\lambda_{i}^{\prime}\|x_{i} -\hat{x}_{i-1}\|^{p+\nu-q}}{\theta_2}\right)^{\frac{\alpha}{1-\alpha}}\gamma,\label{eq:E-9}
\end{eqnarray}
where for $\alpha=1$, we set $\frac{\alpha}{1-\alpha}=\lim_{\alpha\rightarrow 1^{-}}\frac{\alpha}{1-\alpha}=+\infty$.
This would still ensure $\gamma^{\prime}_{i}\in(0, \gamma]$. But notice that $\lambda_{i}^{\prime}$ in the RHS depends on $\gamma^{\prime}_{i}$. To sort out an explicit expression for so-defined $\gamma^{\prime}_{i}$, we denote 
\begin{equation}
\lambda_{i}:= \frac{a_{i}^q}{c_q\gamma A^{q-1}_{i}}. \label{eq:lambda-1}	
\end{equation}
Then by the definition of $\lambda_{i}^{\prime}$ in Lemma \ref{lem:E-3}, with \eqref{eq:E-9} and \eqref{eq:lambda-1}, we can write $\gamma^{\prime}_{i}$ of the form:
\begin{eqnarray}
\gamma^{\prime}_{i} = \left( \frac{L\lambda_{i}  \|x_{i}-\hat{x}_{i-1}\|^{p+\nu-q}}{\theta_2}  \right)^{\alpha}\gamma. \label{eq:gamma-prime-value}
\end{eqnarray}
 
Then by the fact  for all $s\ge 0,t\ge2, x\in\bbR^d$,
\begin{eqnarray}
\|x\|^s\nabla \frac{1}{t}\|x\|^t = \nabla \label{eq:identity} \frac{1}{s+t}\|x\|^{t+s},
\end{eqnarray}
and combing \eqref{eq:lambda-1} and \eqref{eq:gamma-prime-value}, it follows that
\eqref{eq:E-6-2} is equivalent to 
\begin{equation}
\nabla\tilde{f}(x_{i};\hat{x}_{i-1})  + \frac{L^{\alpha}}{q c_q  \lambda_{i}^{(1-\alpha)}\theta_2^{\alpha}} \nabla \frac{1}{\alpha(p+\nu) +(1-\alpha) q}   \|x_{i}-\hat{x}_{i-1}\|^{\alpha(p+\nu) +(1-\alpha) q} = 0.\label{eq:E-10}
\end{equation}
Or equivalently, denote $\varsigma :=\alpha(p+\nu) +(1-\alpha) q$, $x_{i}$ is the solution of the following minimization problem:
\begin{align}
x_{i} := \argmin_{x\in\bbR^d}\left\{ \tilde{f}(x;\hat{x}_{i-1}) +
\frac{L^{\alpha}}{ q c_q  \lambda_{i}^{(1-\alpha)}\theta_2^{\alpha}\varsigma }   \|x-\hat{x}_{i-1}\|^{\varsigma}\right\}.\label{eq:E-7}
\end{align}
In \eqref{eq:E-7}, because $\alpha\in[0,1]$, the power of the norm  $\|x-\hat{x}_{i-1}\|$ ranges from $p+\nu$ to $q$ freely,
which unifies the choice  $\alpha = 0$ in \cite{monteiro2013accelerated} and  $\alpha = 1$ in \cite{bubeck2018near}. 
Meanwhile, \cite{jiang2018optimal, gasnikov2018global} has used a mixture of both $\alpha=0$ and $\alpha=1$ in their formulation, which is also equivalent to \eqref{eq:E-10} by \eqref{eq:identity}.
A surprising phenomenon is that, as our analysis shows, the choice of $\alpha$ in \eqref{eq:E-7} does not affect the convergence rate (except the constant).

 Meanwhile, by  \eqref{eq:gamma-prime-value}, \eqref{eq:E-6-1} is equivalent to 
\begin{eqnarray}
\omega_{i}:=L\lambda_{i}\|x_{i} -\hat{x}_{i-1}\|^{p+\nu-q}\le\theta_2\le1,\label{eq:lambda-2}
\end{eqnarray}
where we call $\omega_{i}$ as a \emph{convergence indicator} in the sense that if for all $1\le i \le k$, $\omega_{i}\le 1$, then the iterate $x_k$ will converge according to the following theorem; otherwise, the convergence of $x_k$ is not guaranteed. 
Then based on the equivalence relationship between \eqref{eq:E-6-2} and \eqref{eq:E-7}, \eqref{eq:E-6-1} and \eqref{eq:lambda-2}, we have Theorem \ref{thm:bound}.
\begin{theorem}[Discrete-Time VPM]\label{thm:bound}
Assume that the convex function $f(x)$ defined in \eqref{eq:prob} satisfies Assumption \ref{ass:f-smooth},  $h(x;x_0)$ satisfies Assumption \ref{ass:h-dis}, $\frac{1}{q}\|\cdot\|^q$ satisfies Assumption \ref{ass:norm}. In \eqref{eq:dis}, $\forall i\ge1,$ the sequences $\{a_i\}, \{A_i\}$ satisfy $a_{i}>0, A_{i} = A_{i-1} +a_{i}$ with $a_0=A_0=0$, $\{x_{i}\}$ satisfies \eqref{eq:E-7},  and $ \{\lambda_{i}\}$ defined in \eqref{eq:lambda-1}	satisfies \eqref{eq:lambda-2}, then for $k\ge 1$, one has
\begin{eqnarray}
f(x_k) - f(x^*)\le \frac{h(x^*;x_0)}{A_k}. \label{eq:bound}
\end{eqnarray}
\end{theorem}
\begin{proof}
See Section \ref{sec:bound}.
\end{proof}
As we see, this theorem is very much like the discrete-time version of the Theorem \ref{thm:cont}. Both try to emulate the convergence rate of the VPM given in Theorem \ref{thm:vpm}. To accurately characterize the convergence rate from \eqref{eq:bound}, we need to have a good lower-bound for $A_k$. However, different from the continuous-time setting, in Theorem \ref{thm:bound}, by $A_{i} = A_{i-1}+a_{i}$, the definition of $\lambda_{i}$ \eqref{eq:lambda-1} and the condition \eqref{eq:lambda-2},  $A_{i}$ must satisfy the condition
\begin{eqnarray}
\frac{L(A_{i}-A_{i-1})^q}{c_q\gamma A^{q-1}_{i}}\|x_{i}-\hat{x}_{i-1}\|^{p+\nu-q}\le \theta_2\le1.\label{eq:A-cond}
\end{eqnarray}
Therefore $A_{i}$ cannot be chosen as an arbitrarily large value as in the continuous-time setting. Except the basic condition $A_0=0$ and for $i\ge 1, A_{i}>0$, \eqref{eq:A-cond} is the only condition $A_{i}$ needs to satisfy, therefore one may expect that the tightest bound of $A_{i}$ should be obtained if $\frac{L(A_{i}-A_{i-1})^q}{c_q\gamma A^{q-1}_{i}}\|x_{i}-\hat{x}_{i-1}\|^{p+\nu-q}=O(1).$ In other words, we hope that 
 \begin{equation}
 0< \theta_1\le L\lambda_{i}\|x_{i}-\hat{x}_{i-1}\|^{p+\nu-q}  = \frac{L(A_{i}-A_{i-1})^q}{c_q\gamma A^{q-1}_{i}}\|x_{i}-\hat{x}_{i-1}\|^{p+\nu-q}\le \theta_2\le 1,
 \end{equation}
where $\theta_1$ and $\theta_2$ are $O(1)$ constants. To verify this point of view, we discuss below the two settings $q = p+\nu$ and $q<p+\nu$, respectively.

When $q = p+\nu$, we have $\lambda_{i} =  \frac{(A_{i}-A_{i-1})^q}{c_q\gamma A^{q-1}_{i}}$ and $\|x_{i}-\hat{x}_{i-1}\|^{p+\nu-q}=1$. Taking $A_{i}$ as a variable, then for all $A_{i} > A_{i-1}$, by the fact $q\ge 2$ and 
\begin{align}
\frac{\d\log \lambda_{i}}{\d A_{i}} = \frac{(q-1)A_{i-1} +A_{i}}{A_{i}(A_{i}-A_{i-1})}> 0,
\end{align}
we have $\lambda_{i}$ is a strictly monotonically increasing function $w.r.t.$ $A_{i}$, which is an one to one mapping. Therefore determining the lower bound of $A_{i}$ is equivalent to determining the lower bound of $\lambda_{i}$. To ensure $E_{i}\le 0$, by the condition \eqref{eq:lambda-2} , when $q = p+\nu$, $L\lambda_{i}$ is upper bounded by the constant $\theta_2\le 1$. Therefore the tightest lower bound for $\lambda_{i}$ is obtained if $L\lambda_{i}$ is lower bounded by a constant $\theta_1\in(0, \theta_2]$. Then by Lemma \ref{lem:seq} and Theorem \ref{thm:bound}, we obtain Theorem \ref{thm:bound-a}.

\begin{theorem}[Convergence Rate for the Case $q=p+\nu$]\label{thm:bound-a}
Assume that the convex function $f(x)$ defined in \eqref{eq:prob} satisfies Assumption \ref{ass:f-smooth},  $h(x;x_0)$ satisfies Assumption \ref{ass:h-dis}, $\frac{1}{q}\|\cdot\|^q$ satisfies Assumption \ref{ass:norm}. $c_q$ is defined in Lemma \ref{lem:E-3}. In \eqref{eq:dis}, $\forall i\ge1,$ the sequences $\{a_i\}, \{A_i\}$ satisfy $ a_{i}>0, A_{i} = A_{i-1} +a_{i}$ with $a_0= A_0 = 0$, $\{x_{i}\}$ satisfies \eqref{eq:E-7},  and $ \{\lambda_{i}\}$ defined in \eqref{eq:lambda-1} satisfies
\begin{eqnarray}
0<\theta_1\le L\lambda_{i}\le \theta_2\le 1, \label{eq:lambda-cond}
\end{eqnarray}
then for $k\ge 1$, we have
\begin{eqnarray}
A_k\ge \frac{\theta_1 c_q \gamma}{L} \left(\frac{k}{p+\nu}\right)^{p+\nu},
\end{eqnarray}
and
\begin{eqnarray}
f(x_k) - f(x^*)\le  \frac{h(x^*;x_0)}{A_k}\le   \frac{L}{\theta_1 c_q \gamma}h(x^*;x_0) \left(\frac{p+\nu}{k}\right)^{p+\nu}. \label{eq:bound-a}
\end{eqnarray}
\end{theorem}
\begin{proof}
See Section \ref{sub:thm:bound-a}.
\end{proof}

\medskip

When $q< p+\nu$, because the  condition of $\lambda_{i}$ to ensure $E_{i}\le 0$ involves the unknown $x_{i}$, the situation seems to be more complicated. 
Nevertheless, under the conditions \eqref{eq:E-7} and \eqref{eq:lambda-2}, and combining Lemmas \ref{lem:dis-upper} and \ref{lem:dis-lower}, we can obtain a condition as in Lemma \ref{lem:q-neq-p-nu} below that leads to a good lower bound for $A_k$.
\begin{lemma}\label{lem:q-neq-p-nu}
Assume  $\{x_{i}\}$ satisfies \eqref{eq:E-7}, $\{\omega_{i}\}$ satisfies \eqref{eq:lambda-2}. Then if $2\le q < p+\nu$, we have
\begin{align}
 \sum_{i=1}^{k}\omega_{i}^{\frac{\varsigma}{p+\nu-q}}
\left(\frac{A_{i}^{p+\nu-1}}{ (A_{i} -A_{i-1})^{p+\nu}}\right)^{\frac{q}{p+\nu-q}}\!\!\!\le
{q\theta_2^{\alpha}} (1-\theta_2^{\frac{q}{q-1}})^{-1}\gamma^{-\frac{p+\nu}{p+\nu-q}}\left(\frac{L}{c_q}\right)^{\frac{q}{p+\nu-q}}h(x^*;x_0).\label{eq:q-neq-p-nu-4}
\end{align}
\end{lemma}
\begin{proof}
See Section \ref{sub:lem:q-neq-p-nu}.
\end{proof}

In Lemma \ref{lem:q-neq-p-nu}, if $\theta_2\in(0,1)$, then the RHS of \eqref{eq:q-neq-p-nu-4} will be a positive constant. Therefore  \eqref{eq:q-neq-p-nu-4} will have the same form of \eqref{eq:seq-3} of Lemma \ref{lem:seq} if $\omega_{i}$ on the LHS of \eqref{eq:q-neq-p-nu-4} is lower bounded by a constant $\theta_1\in(0, \theta_2].$ Based on the above analysis, and combining  Lemma \ref{lem:seq}, Theorem \ref{thm:bound} and Lemma \ref{lem:q-neq-p-nu}, we can characterize the convergence rate of the proposed iteration when $2\le q<p+\nu$.  
\begin{theorem}[Convergence Rate for the Case $2\le q<p+\nu$]\label{thm:bound-b}
Assume that the convex function $f(x)$ defined in \eqref{eq:prob} satisfies Assumption \ref{ass:f-smooth}, and  $h(x;x_0)$ satisfies Assumption \ref{ass:h-dis}, $\frac{1}{q}\|\cdot\|^q$ satisfies Assumption \ref{ass:norm}.  $c_q$ is defined in Lemma \ref{lem:E-3}. In \eqref{eq:dis}, $\forall i\ge1,$ the sequences $\{a_i\}, \{A_i\}$ satisfy $ a_{i}>0, A_{i} = A_{i-1} +a_{i}$ with $a_0=A_0=0$, $\{x_{i}\}$ satisfies \eqref{eq:E-7},  and $ \{\lambda_{i}\}$ defined in \eqref{eq:lambda-1} satisfies
\begin{eqnarray}
0<\theta_1\le \omega_{i}=L\lambda_{i}\|x_{i}-\hat{x}_{i-1}\|^{p+\nu-q}\le \theta_2<1, \label{eq:lambda-cond-a}
\end{eqnarray}
then define $C_0:=\left({q\theta_2^{\alpha}} \big(1-\theta_2^{\frac{q}{q-1}}\big)^{-1}\right)^{-\frac{p+\nu-q}{q}} \theta_1^{\frac{\alpha(p+\nu) + (1-\alpha)q}{q}}\gamma^{\frac{p+\nu}{q}}{c_q}$,  we have

\begin{eqnarray}
A_k \ge \frac{C_0}{L}(h(x^*;x_0))^{-\frac{p+\nu-q}{q}}\left(\frac{k}{p+\nu}\right)^{\frac{(q+1)(p+\nu)-q}{q}}
\end{eqnarray}
and
\begin{eqnarray}
f(x_k) - f(x^*)\le  \frac{h(x^*;x_0)}{A_k}\le \frac{L}{C_0}(h(x^*;x_0))^{\frac{p+\nu}{q}} \left(\frac{p+\nu}{k}\right)^{\frac{(q+1)(p+\nu)-q}{q}}.     \label{eq:bound-b}
\end{eqnarray}
\end{theorem}

\begin{proof}
See Section \ref{sub:thm:bound-b}.
\end{proof}

In Theorems \ref{thm:bound-a} and \ref{thm:bound-b}, if we do not consider the constants, in both $q=p+\nu$ and $2\le q<p+\nu$ settings, we can find an $\epsilon$-accurate solution $x$ such that $f(x)- f(x^*)\le \epsilon$ with  $$O\Big(\epsilon^{-\frac{q}{(q+1)(p+\nu)-q}}\Big)$$ iterations, where $q\in[2, p+\nu]$. It is easy to find that the rate will be the best as $O\big(\epsilon^{-\frac{2}{3(p+\nu)-2}}\big)$ if we set $q=2$. In fact $O\big(\epsilon^{-\frac{2}{3(p+\nu)-2}}\big)$ matches the lower bound of iteration complexity \cite{grapiglia2019tensor}
 for all the settings of $p\in\{1,2,\ldots, \}$ and $\nu\in[0,1]$ as long as $p+\nu\ge2$. As $q$ becomes large, the rate $O\Big(\epsilon^{-\frac{q}{(q+1)(p+\nu)-q}}\Big)$ will become worse. However, particularly, when $q = p+\nu$, $\lambda_{i}$ can be determined trivially and thus the setting $q = p+\nu$ is suboptimal but has the advantage of algorithmic implementation, as we will elaborate on later.

Regarding the other two parameters $\theta_1, \theta_2$, when $q=p+\nu$, based on Theorem \ref{thm:bound-a},   to minimize the bound in \eqref{eq:bound-a}, the optimal choice will be $\theta_1=1$ and thus $\theta_2 =1$ by $\theta_1\le \theta_2\le1$. When $q<p+\nu$, based on Theorem \ref{thm:bound-b}, one can optimize the choice of $\theta_1, \theta_2$ by minimizing the bound in \eqref{eq:bound-b} under the constraint $0<\theta_1\le \theta_2<1$. However, for the the case with $2\le q<p+\nu$, the choice of $\theta_1, \theta_2$ also influences the complexity to find $\lambda_{i}$ that satisfies \eqref{eq:lambda-cond-a}.

As we have noted before, by varying the parameter $\alpha$ from $0$ to $1$ in \eqref{eq:E-7}, the range of the power of $\|x-\hat{x}_{i-1}\|$ changes from $q$ to $p+\nu$. For $ q=p+\nu$, as Theorem \ref{thm:bound-a} indicates, choice of $\alpha$ has no influence on the convergence rate; for $2\le q<p+\nu$, as Theorem \ref{thm:bound-b} shows, $\alpha$ only has a minor influence on the constant in the bound. Therefore, our result shows that $\alpha$ can be chosen according to implementation convenience without worrying about the convergence rate.

Compared with the existing papers about high-order optimization \cite{nesterov2006cubic,nesterov2008accelerating,nesterov2018implementable,grapiglia2019tensor} and \cite{monteiro2013accelerated,bubeck2018near,jiang2018optimal,gasnikov2018global}, our convergence results are given under the H\"{o}lder continuous assumption $w.r.t.$ a general norm $\|\cdot\|$ that satisfies Assumption \ref{ass:norm}. Such general norms include the Euclidean norm $\|x\|_2$ and the generalized Euclidean norm $\sqrt{x^TBx}$ as special cases, where $B$ is any positive definite matrix. To the best of our knowledge, this is the first convergence result for high-order optimization that can be applied to the high-order non-Euclidean smoothness setting. To this end, we have adopted a new proof paradigm inspired by the intuitive proof techniques of AXGD \cite{diakonikolas2017accelerated} for first order methods.

Summarizing the above results, we obtain a \emph{unified acceleration framework} (UAF) shown in Algorithm \ref{alg:dis}. In Algorithm \ref{alg:dis},  the parameters $p, \nu$ are from the problem setting and the parameters $q, \alpha$ and the proxy function $h(x;x_0)$ are for framework design. These parameters can vary in their entire feasible ranges. By specifying $p, \nu, q, \alpha$ and $h(x;x_0)$, we obtain algorithmic instances of UAF. 
As results, Algorithm \ref{alg:dis} recovers many existing algorithms. We give a few examples in Table \ref{tab:unified}\footnote{As we have mentioned, \cite{jiang2018optimal,gasnikov2018global} have used a mixture of $0$ and $1$ for $\alpha$, which is also equivalent to  \eqref{eq:E-10} by \eqref{eq:identity}. To simplify the presentation and also the consideration of practical implementation, we only consider the choice of a single $\alpha$ in Algorithm \ref{alg:dis}. }.

Meanwhile, Algorithm \ref{alg:dis} also includes several new interesting instances. First, if we set $p=\nu=1, q=2,\alpha\in[0,1]$, Algorithm \ref{alg:dis} defines a new variant of AGD with an $O(1/k^2)$ convergence rate. Such variant is similar to a variant accelerated extra-gradient descent (AXGD) of AGD. One advantage of this variant is that Algorithm \ref{alg:dis} allows $h(x;x_0)$ to be any strongly convex function $w.r.t.$  $\|\cdot\|$, while AXGD assumes that $h(x;x_0)$ is the Bregman divergence of a strongly convex function $w.r.t.$ $\|\cdot\|$. Second, if we set $p\in\{2,3,\ldots, \}, \nu\in[0,1], q=2, \alpha\in[0,1]$, then we obtain the first kind of high-order algorithms that can attain the optimal rate $O(\epsilon^{-\frac{2}{3(p+\nu)-2}})$ for the composite minimization problem \eqref{eq:prob} with the smooth component $g(x)$ having $(p, \nu, L)$-H\"{o}lder continuous derivatives $w.r.t.$ $\|\cdot\|.$

\begin{algorithm}[t!]
\caption{Unified Acceleration Framework (UAF)}\label{alg:dis}
\begin{algorithmic}[1]
\STATE \textbf{Input:} a convex function $f(x)$ defined in \eqref{eq:prob}  satisfies Asumption \ref{ass:f-smooth}; $\hat{f}(x;y)$ and  $\tilde{f}(x;y)$ defined in Lemma \ref{lem:hat-tilde-f}.
$h(x;x_0)$ satisfies Assumption \ref{ass:h-dis};
$\frac{1}{q}\|\cdot\|^q$ satisfies Assumption \ref{ass:norm}.

\STATE Set  constants $\theta_1, \theta_2$ such that  $\theta_1\in(0, \theta_2]$, $\theta_2\in(0,1]$ if $q=p+\nu$, $\theta_2\in(0,1)$ if $2\le q<p+\nu$.
\STATE Set $\alpha\in[0,1], c_q=\left( {\beta(q-1)^{1-q}} \right)^{\frac{1}{q}}, \varsigma = \alpha(p+\nu) +(1-\alpha) q.$
\STATE Set $A_0 = 0, x_0 = z_0\in\bbR^d$. 

\FOR{ $i = 1$ \TO $ k $ }

\STATE $x_{i} = \argmin_{x\in\bbR^d}\left\{ \tilde{f}(x;\hat{x}_{i-1}) +   \frac{L^{\alpha}}{ q c_q  \lambda_{i}^{(1-\alpha)}\theta_2^{\alpha}\varsigma }    \|x-\hat{x}_{i-1}\|^{\varsigma}\right\}$,

where we find a $\lambda_{i}>0$ such that $a_{i}, A_{i}, \lambda_{i}$  and $\hat{x}_{i-1}\in\bbR^d$ satisfying
\begin{eqnarray}
&& A_{i} = A_{i-1} + a_{i},\; \lambda_{i}= \frac{a_{i}^{q}}{c_q \gamma A_{i}^{q-1}},\; \hat{x}_{i-1} = \frac{A_{i-1}}{A_{i}}x_{i-1} + \frac{a_{i}}{A_{i}}z_{i-1},\label{eq:lam-A-x}\\
 && \theta_1 
\le L \lambda_{i}\|x_{i} - \hat{x}_{i-1}\|^{p+\nu-q}\le \theta_2.\label{eq:q-equal-p-nu}
\end{eqnarray}

\STATE Update $z_{i} =  \argmin_{x\in \bbR^d}\left\{\sum_{j=1}^{i}a_{j}\hat{f}(x; x_{j})+h(x; x_0)
\right\}$.
\ENDFOR
\RETURN $x_k$.
	\end{algorithmic}	
\end{algorithm}

\begin{table}[t]
{\footnotesize
  \caption{Algorithmic Instances of the Unified Acceleration Framework with $h(x;x_0):=\frac{1}{q}\|x-x_0\|_2^q$.}\label{tab:unified}
\begin{center}
  \begin{tabular}{|c|c|c|c|c|} \hline
   Instances &  $p$ & $\nu$ & $q$ &  $\alpha$   \\ \hline
    \cite{baes2009estimate,nesterov2018implementable} & $\{2,3,\ldots\}$ &  $1$ & $p+1$ & $1$ \\
     \cite{grapiglia2019tensor}              &   $\{2,3,\ldots\}$ & $[0,1]$ & $p+\nu$ & $1$ \\  
    \cite{monteiro2013accelerated}  &  $2$  & $1$ & $2$ & $0$   \\ 
    \cite{bubeck2018near}                   &  $\{2,3,\ldots\}$ & $1$ & $2$ & $1$ \\ 
    \cite{jiang2018optimal,gasnikov2018global} & $\{2,3,\ldots\}$ & $1$ & $2$ & a mixture of $0$ and $1$ \\  
    \hline
  \end{tabular}
\end{center}
}
\end{table}

Note that in Algorithm \ref{alg:dis}, under the Assumptions \ref{ass:h-dis}, \ref{ass:norm} and \ref{ass:f-smooth}, in Step 2, we choose $3$ design parameters $\theta_1, \theta_2, \alpha$ to be used in the Algorithm.  For the loop from Step 4 to 7, there are mainly two subproblems to solve:
\begin{itemize}
    \item The first one is about finding  $\lambda_{i}$ such that the minimizer $x_{i}$ of the objective \eqref{eq:E-7}, together with $\lambda_{i}$, satisfy the conditions \eqref{eq:lam-A-x} and \eqref{eq:q-equal-p-nu}. 
    \item The second one is about finding the solution $z_{i}$ of a discrete-time convex approximation problem of the VPM in Step $6$. Because in our setting the convex approximation $\hat{f}(x;y)$ defined in Lemma \ref{lem:hat-tilde-f} is a linear function plus a simple convex function $l(x)$, the subproblem of finding $z_{i}$ can be solved efficiently. 
\end{itemize} 

When $p=\nu=1$ and $q=2$, the subproblem associated with Step 5, namely \eqref{eq:E-7}, is reduced to a proximal gradient decent step \cite{parikh2014proximal}, which can be solved efficiently. However, in the setting of high-order optimization, $i.e., p\in\{2, 3,\ldots\}$, \eqref{eq:E-7} is nontrivial and will dominate the per-iteration cost in general. Finding a general efficient procedure to solve this subproblem remains active research. Nevertheless, for some special important cases, there already exist efficient algorithms. For example, if $p=2, \nu=1, \alpha = 1$ and the maybe nonsmooth part $l(x)=0$, \eqref{eq:E-7} is reduced to an iteration of cubic regularized Newton method (CNM), which can be solved efficiently by the Lanzcos method \cite{carmon2018analysis}; if $p=3, \nu=1, \alpha =1$ and $l(x)=0$, \eqref{eq:E-7} is reduced to a third-order convex multivariate polynomial and can be solved as efficiently as the iteration of CNM in many cases \cite{nesterov2018implementable,bullins2018fast}.

Notice that, in Step 5, for the setting $q=p+\nu$, $\lambda_{i}$ can be determined easily as it does not depend on $x_{i}$ and thus $A_{i}, a_{i}$ can be solved efficiently by solving a simple one-dimensional equation with Newton method. However,
for the setting $2\le q<p+\nu$, the condition \eqref{eq:q-equal-p-nu} depends on the solution $x_{i}$ and can not be determined so trivially. In fact, as of now, when $2\le q<p+\nu$, we do not even know whether we can find such a pair $(x_{i}, \lambda_{i})$ that satisfies all the conditions simultaneously.
As nearly a trivial extension to \cite{bubeck2018near}, the following Proposition \ref{prop:exist} ensures  such a pair always exists until we attain the minimizer.
\begin{proposition}\label{prop:exist}
Let $A\ge 0, \lambda\ge 0, x,y\in\bbR^d$ such that $f(x)\neq f(x^*)$. Assume that $a(\lambda)$ is implicitly defined by  
\begin{eqnarray}
\lambda = \frac{(a(\lambda))^q}{c_q \gamma (A+a(\lambda))^{q-1}},\quad\text{\rm and}\quad x(\lambda) = \frac{a(\lambda)}{A+a(\lambda)} x+\frac{A}{A+a(\lambda)} y,
\end{eqnarray}
\begin{align}
w(\upsilon) &= \argmin_{z\in\bbR^d}\left\{ \tilde{f}(z;\upsilon) +   \frac{L^{\alpha}}{ q c_q  \lambda^{(1-\alpha)}\theta_2^{\alpha}\varsigma}   \|z-\upsilon\|^{\varsigma}\right\}\label{eq:y-z},\\
\chi(\lambda) &= L\lambda\|w(x(\lambda))-x(\lambda)\|^{p+\nu -q},  
\end{align}
where the constants $p, q,\nu, \alpha,c_q, \gamma, L, \varsigma$ and $\theta_2$ are given in Algorithm  \ref{alg:dis}. Then $\chi(\lambda)$ is a continuous function with $\chi(0) =0 $ and $\chi(+\infty)=+\infty.$
\end{proposition}
\begin{proof}
See Section \ref{sec:prop:exist}
\end{proof}
By Proposition \ref{prop:exist}, with the setting $A:=A_{i-1}, x := z_{i-1}, y:=x_{i-1}$, we can always use a binary search procedure to find a pair $(x_{i}, \lambda_{i})$ such that $\chi(\lambda_{i}) = L\lambda_{i}\|x_{i} - \hat{x}_{i-1}\|^{p+\nu-q}$ satisfies the condition \eqref{eq:q-equal-p-nu}. For the case with $\alpha =0, q=2$ and $\|\cdot\|:=\|\cdot\|_2$, a complexity analysis for a binary search procedure can be found in \cite{jiang2018optimal}; for the case with $\alpha=1, p\in\{2,3,\ldots\}, \nu=1$ and $\|\cdot\|:=\|\cdot\|_2$, a   complexity analysis for a binary search procedure can be found in \cite{bubeck2018near}. Although it is possible to give a complexity analysis of binary search for the general setting in \eqref{eq:y-z}, in this paper we consider another perspective. 

In the Discussion section of \cite{nesterov2018implementable}, Nesterov claims that from the view of practical efficiency, the Algorithm  \ref{alg:dis} with the suboptimal setting $q = p+\nu$ may be better than the Algorithm  \ref{alg:dis} with the optimal setting $q=2$, where ``optimal'' is in the sense of iteration complexity. If we do not consider the implementation cost in the Step 5 of Algorithm  \ref{alg:dis} and ignore the difference of constants in the bound of Theorems \ref{thm:bound-a} and \ref{thm:bound-b}, to attain an $\epsilon$-accurate solution such that $f(x)-f(x^*)\le \epsilon$, the ratio from the number of iterations of the suboptimal algorithm with $q=p+\nu$ to that of the optimal algorithm with $q=2$ is $$O\left(\left(\frac{1}{\epsilon}\right)^{\frac{1}{p+\nu}-\frac{2}{3(p+\nu)-2)}}\right)=O\left(\left(\frac{1}{\epsilon}\right)^{\frac{p+\nu-2}{(p+\nu)(3(p+\nu)-2)}}\right).$$ 
If $p=2, \nu=1$, $i.e.,$ the commonly second-order setting, the ratio will be $O\left( \left(\frac{1}{\epsilon}\right)^{\frac{1}{21}} \right)$, which implies that when we pursue an accuracy $\epsilon = 2^{-21}\approx 10^{-6}$, if the per-iteration cost of the optimal setting $q = 2$ (or the settings $2\le q< p+\nu$) is twice larger than the suboptimal setting $q = p+\nu$, then the small advantage of the optimal setting will be removed by the additional implementation complexity. Because of this effect, a binary search procedure which involves $O(\log \frac{1}{\epsilon})$ calls to the subprocedure of finding $x_{i}$ may be rather unrealistic in practice. Therefore in this paper, instead of binary search, we propose a simple heuristic to find a pair $(x_{i}, \lambda_{i})$ to satisfy the condition \eqref{eq:q-equal-p-nu}. The proposed heuristic  
only needs $1$ call to the subprocedure of finding $x_{i}$ and will be clear in Section \ref{sec:exp}.

\begin{remark}
 The idea that two-step fixed-point iterations lead to acceleration is first introduced in \cite{diakonikolas2017accelerated}, which has proposed a variant accelerated extra-gradient descent (AXGD) of AGD. In this paper, such point of view motivates us to simplify the analysis by defining an intermediate variable $\hat{z}_{i}$ in \eqref{eq:two-fixed}, while the strategy leading to acceleration in this paper is using a combination of a convex approximation \eqref{eq:dis} of the original VPM problem \eqref{eq:vpm} and a smooth approximation \eqref{eq:E-7} of the intermediate VPM
 problem \eqref{eq:E-3}.
\end{remark}

\medskip

\begin{remark}
Similar to \cite{grapiglia2019tensor}, it is also possible to give a universal version of Algorithm \ref{alg:dis} in the sense that following the paradigm of \cite{grapiglia2019tensor}, with corresponding modifications of Algorithm \ref{alg:dis}, we can obtain a near-optimal rate even if the H\"{o}lder parameter $\nu$ is unknown. Such modification is of interest, however it goes beyond the topic of this paper and will be left for further research.
\end{remark}

\section{General Restart Scheme for Uniformly Convex Functions}\label{sec:restart}
Algorithm \ref{alg:dis} is  proposed mainly for minimizing convex functions, which matches the lower bounds of iteration complexity given in \cite{arjevani2017oracle,grapiglia2019tensor} for the class of functions considered. Nevertheless, for uniformly convex functions, we should expect better convergence rates, as shown in \cite{arjevani2017oracle}. It is however nontrivial to  match such lower bounds by following the techniques introduced in the previous section. In this section, inspired by the {\em restart method} for accelerated cubic regularized Newton method (ACNM) \cite{nesterov2008accelerating}, we propose a  \emph{general restart scheme} for any general algorithm $\cA$ with a specified form of convergence rate. Then we apply the proposed restart scheme to the Algorithm \ref{alg:dis} and obtain convergence rates that match the lower bound given in \cite{arjevani2017oracle}.

To describe the restart scheme, we first make some assumptions about the function and the algorithm of consideration.

\begin{assumption}\label{ass:s-uniform}
  $f(x)$  is $(s, \sigma)$-uniformly convex $w.r.t.$ $\|\cdot\|$, where $s\ge 2, \sigma >0$.
\end{assumption}

\begin{assumption}\label{ass:alg}
To minimize $f(x)$, let $\cA_m(y)$ denote the output of an algorithm $\cA$ after $m$ iterations with an input $y\in\bbR^d$. We assume that the output satisfies 
\begin{eqnarray}
f(\cA_m(y)) - f(x^*)\le \frac{c_{\cA }\|y-x^*\|^{v}}{m^r}, \label{alg:form} 
\end{eqnarray}
for some constants $r>0, v>0, c_{\cA}>0$. 
\end{assumption}

Let $R> 0$ be a constant such that $\|x_0-x^*\|\le R$, for the initial point $x_0\in \bbR^d$ and $x^* \in \bbR^d$ a minimum point of $f(x)$.  We define two constants: 
\begin{equation}
 m_0 = \left\lceil \left(\frac{2^s s c_{\cA} R^{v-s}}{\sigma}\right)^{\frac{1}{r}}\right\rceil, \quad 
 k_0=
 \begin{cases}
 \left\lceil \frac{1}{s} + \frac{v}{s}\log_2 R+   \frac{1}{v-s}\log_2 \left( \frac{s c_{\cA}}{\sigma}\right)\right\rceil\,  &\mbox{if} \, s< v, \\
 +\infty  &\mbox{if} \, s\ge v
 \end{cases}.
\label{eq:m-0-k}
\end{equation}
 Here $m_0$ and $k_0$ are carefully chosen for consideration of best convergence rates (as one will see in the proof of the theorem about the convergence rates).

Then given a total number of epoches $K \in \mathbb{Z}_+$, Algorithm \ref{alg:restart} gives a general restart scheme  for minimizing uniformly convex functions.
\begin{algorithm} [ht!]
\caption{A General Restart Scheme}\label{alg:restart}
\begin{algorithmic}[1]
\STATE {\bf Input:} An $(s, \sigma)$-uniform convex function $f(x)$ with $ s\ge 2, \sigma>0$, an algorithm $\cA$ satisfying \eqref{alg:form}, and a total number of epoches $K \in \mathbb{Z}_+$.

\STATE {\bf Set} $y_0 = x_0 \in \mathbb{R}^d$ and $m_0, k_0$ defined in \eqref{eq:m-0-k}.
\FOR{$k = 0,1,2,\ldots,K-1$}
\IF{$k\le k_0-1$}
\STATE $m_k = \left\lceil m_{0} 2^{-\frac{v-s}{r} k}\right\rceil$.
\STATE $y_{k+1} = \mathcal{A}_{m_k}(y_k)$.
\ELSE
\STATE $y_{k+1}=\mathcal{A}_1(y_k).$
\ENDIF
\ENDFOR
\RETURN $y_K$.
	\end{algorithmic}	
\end{algorithm}

As we see in the algorithm, in each epoch, from Step 4 to 6, we set the number of iterations to be $m_k$ and update the iterate $y_k$ by calling the inner algorithm $\cA$ with $m_k$ iterations. 
If $s=v$, then $m_k$ will remain as a constant along the epochs; if $s>v$, then $m_k$ will increase by a linear rate. In the settings $s\ge v$, $k_0=+\infty$ and therefore the steps from 7 to 9 in Algorithm \ref{alg:restart} will not be executed. If $s<v$, $m_k$ will decrease by a linear rate until $k=\min\{K,k_0\}-1$. After $k>k_0$, in each epoch the number of iterations of calling $\cA$ is set to be $1$. With these settings, the convergence behavior of Algorithm \ref{alg:restart} in terms of the number of epoches $k$ is given by the following theorem. 
\begin{theorem}[Convergence Rates with the General Restart Scheme]
\label{thm:restart}
In Algorithm \ref{alg:restart}, when $k\le k_0$ , it follows that
\begin{eqnarray}
f(y_k) - f(x^*)\le \frac{\sigma}{s 2^{sk}}\|x_0-x^*\|^{v}.\label{eq:K-lower}
\end{eqnarray}
If $s<v$ and when $k_0<k\le K$, one has
\begin{eqnarray}
f(y_{k}) - f(x^*)\le \left(\frac{\sigma^{v}}{s^v c_{\cA}^s}\right)^{\frac{1}{v-s}} 2^{-(\frac{v}{s})^{k-k_0}}. \label{eq:K-upper}
\end{eqnarray}
\end{theorem}
\begin{proof}
See Section \ref{sec:thm:restart}.
\end{proof}
By Theorem \ref{thm:restart}, in the first stage when $k\le k_0$, Algorithm \ref{alg:restart} will converge at a linear rate in terms of the number of epochs. 
If $s<v$, in the second stage when $k>k_0$, it will converge fast at a superlinear rate $O(2^{-(\frac{v}{s})^{k-k_0}})$.

In the first stage, the number of iterations $m_k$ of calling $\cA$ will depend on the relation between the two parameters $s$ and $r$. The convergence rate in terms of the total number of iterations of calling $\cA$ will be very different for the three settings $s=v$, $s<v$ and $s>v$. If $s=v$, for $k\le k_0-1,$ the total iterations will be $\sum_{i=1}^k m_i = k m_0$, which is increasing as a linear function of the epoch $k$; if $s<v$, it will increase exponentially in $k$; if $s>v$, $\sum_{i=1}^k m_i$ will remain as a constant if we do not consider the ceiling function. 

Through easy computation, we have the following corollary for measuring the complexity of Algorithm \ref{alg:restart} in terms of the total iterations of calling the algorithm $\cA$. 
\begin{corollary}\label{cor:1}
In Algorithm \ref{alg:restart}, to find an $\epsilon$-accurate solution such that $f(x)-f(x^*)\le \epsilon$,  if $s=v$, the total number of iterations of calling $\cA$  is at most 
\begin{eqnarray}
O\left(\left( \frac{c_{\cA}}{\sigma}\right)^{\frac{1}{r}}\log_2 \frac{\sigma R^v}{\epsilon}\right). \label{eq:linear-res-A-2}
\end{eqnarray} 
If $s<v$, then the total number of iterations of calling $\cA$ is at most
\begin{eqnarray}
O\left(\left(\frac{c_{\cA R^{v-s}}}{\sigma} \right)^{\frac{1}{r}}+\log_{\frac{p+\nu}{s}}\log_{2}\left(\left(\frac{\sigma^{v}}{c_{\cA}^s} \right)^{\frac{1}{v-s}}\frac{1}{\epsilon}\right)\right).  \label{eq:super-res-A}
\end{eqnarray}
If $s>v$, the total number of iterations of calling $\cA$  is at most 
\begin{eqnarray}
O\left(c_{\cA}^{\frac{1}{r}}R^{-\frac{(s-v)^2}{sr}} \sigma^{-\frac{v}{sr}}\epsilon^{-\frac{s-v}{sr}}\right).\label{eq:sublinear-res-A}
\end{eqnarray}
\end{corollary}

In the above discussion, we only assume that $\cA$ in Assumption \ref{ass:alg} has a convergence rate of the form \eqref{alg:form}, while  $c_{\cA}>0, r>0, v>0$ are unspecified parameters. Such a form of convergence rates appears widely for both first-order algorithms \cite{nesterov1998introductory} and high-order algorithms (as shown in Theorems \ref{thm:bound-a} and \ref{thm:bound-b} if we set $h(x;x_0):=\frac{1}{q}\|x-x_0\|^q$). 
Therefore, although the above restart framework is mainly proposed to restart our Algorithm \ref{alg:dis} for high-order optimization, it may be of independent interest for other algorithms. Below are some examples of Algorithm \ref{alg:restart} in the first-order settings.

\begin{example}
In the nonsmooth setting, let $\cA$ denote gradient descent methods \cite{nesterov1998introductory} or \cite{hazan2016introduction}. Then with suitable parameter settings $\cA$ has a rate of the form $O\big(\frac{c_{\cA}\|x_0-x^*\|}{\sqrt{m}}\big)$.
If the function $f(x)$ to be minimized satisfies Assumption \ref{ass:s-uniform}, then by restarting $\cA$ with Algorithm \ref{alg:restart} and according to \eqref{eq:sublinear-res}, we obtain a rate $O\Big({\sigma}^{-\frac{2}{s}}\epsilon^{-\frac{2(s-1)}{s}}\Big)$, which matches the lower bound of nonsmooth and uniformly convex optimization \cite{juditsky2014deterministic}\footnote{In fact, the multi-stage algorithm in \cite{juditsky2014deterministic} is a special case of Algorithm \ref{alg:restart} with the setting $v=1, r = \frac{1}{2}$}.  Particularly, the obtained rate recovers the well-known lower bound $O(\sigma^{-1}\epsilon^{-1})$ of nonsmooth and strongly convex optimization \cite{juditsky2014deterministic} by setting $s=2$. 
\end{example}

\begin{example}
In the smooth setting, let $\cA$ denote accelerated gradient descent methods  \cite{nesterov1998introductory}, then $\cA$ has a rate of the form $O\Big(\frac{L\|x_0-x^*\|^2}{m^2}\Big)$, where $L$ denotes the smoothness constant. If the function $f(x)$ to be minimized satisfies Assumption \ref{ass:s-uniform} with $s=2$ ($i.e.,$ strongly convex), then by restarting $\cA$ with Algorithm \ref{alg:restart} and according to \eqref{eq:linear-res-A-2}, we obtain a rate $O\Big(\left(\frac{L}{\sigma}\right)^{\frac{1}{2}}\log_2 {\frac{1}{\epsilon}} \Big)$, which matches the lower bound of smooth and strongly convex optimization \cite{nesterov1998introductory}. 
\end{example}

Now, let us consider how to further improve the proposed UAF Algorithm \ref{alg:dis} in the uniformly convex setting. By applying the proposed restart Algorithm \ref{alg:restart}, we expect to obtain better convergence rates. The following lemma and theorem below make this precise.

\begin{lemma}\label{lem:simplified-bound}
When $\cA$  in Assumption \ref{ass:alg} is Algorithm \ref{alg:dis} with its parameter settings and $h(x;x_0):=\frac{1}{q}\|x-x_0\|^q$, then we have
\begin{eqnarray}
f(\cA_m(y)) - f(x^*) \le {c_{\cA}\| y-x^*\|^{p+\nu}}{m^{-\frac{(q+1)(p+\nu)-q}{q}}},
\end{eqnarray}
where
\begin{equation}
c_{\cA}=
\begin{cases}
\frac{1}{\theta_1 c_q \gamma} \left({p+\nu}\right)^{p+\nu}L & \;{\rm if}\;\; q = p+\nu, \\
C_0^{-\frac{p+\nu-q}{q}} \left({p+\nu}\right)^{\frac{(q+1)(p+\nu)-q}{q}}L & \;{\rm if}\;\; 2\le q < p+\nu,
\end{cases}
\end{equation}
where the constants such as $ \theta_1,\theta_2, c_q, L, C_0$ are from Theorems \ref{thm:bound-a} and \ref{thm:bound-b}.
\end{lemma}

Then by Lemma \ref{lem:simplified-bound}, if $f(x)$ satisfies Assumption \ref{ass:s-uniform}, then we can use the restart framework in Algorithm \ref{alg:restart} to further accelerate Algorithm \ref{alg:dis} to obtain better rates in the uniformly convex setting. In fact, by setting $r = \frac{(q+1)(p+\nu)-q}{q}$ in Corollary \ref{cor:1}, we directly obtain the following result.
\begin{theorem}[UAF with Restart]
For a $(s, \sigma)$-uniformly convex and $(p, \nu, L)$-H\"{o}lder continuous function, with the proposed restart scheme, Algorithm \ref{alg:restart} applied to the unified acceleration framework Algorithm \ref{alg:dis}, to find an $\epsilon$-accurate solution, the number of iterations we need is at most
\begin{eqnarray}
O\left(\left( \frac{L}{\sigma}\right)^{\frac{q}{(q+1)(p+\nu)-q}}\log \frac{1}{\epsilon}\right) \; \text{if} \; s = p + \nu; \label{eq:linear-res-p2}
\end{eqnarray} 
\begin{eqnarray}
O\left(\left(\frac{L}{\sigma} \right)^{\frac{q}{(q+1)(p+\nu)-q}}+\log\log\left(\left(\frac{\sigma^{p+\nu}}{L^s} \right)^{\frac{1}{p+\nu-s}}\frac{1}{\epsilon}\right)\right) \; \text{if} \; s < p + \nu; \label{eq:super-res-p2}
\end{eqnarray}
\begin{eqnarray}
O\left(\left(  
\frac{L}{\sigma}\right)^{\frac{q}{(q+1)(p+\nu)-q}}\left(\frac{\sigma}{\epsilon}\right)^{\frac{(s-p-\nu)q}{s((q+1)(p+\nu)-q)}}\right) \; \text{if} \; s > p + \nu..
\label{eq:sublinear-res-p2}
\end{eqnarray}
\end{theorem}
These were the rates  \eqref{eq:linear-res}, \eqref{eq:super-res}, and \eqref{eq:sublinear-res} given earlier in the Introduction Section \ref{sec:intro}.

\section{Implementation Details and Experimental Validation}\label{sec:exp}
In high-order optimization, a very different situation from first-order optimization is that the optimal acceleration method ($e.g.,$ the UAF with $q=2$) requires certain conditions (in each iteration). Those are not so  trivial to compute or be satisfied. 
In fact, in the UAF Algorithm \ref{alg:dis}, for $2\le q<p+\nu$, it is not trivial to find a pair $(\lambda_{i}, x_{i})$ to satisfy the condition \eqref{eq:q-equal-p-nu}. 
In \eqref{eq:lambda-2}, we have defined $\omega_{i}=L \lambda_{i}\|x_{i} - \hat{x}_{i-1}\|^{p+\nu-q}$ as a convergence indicator in the sense that  $\forall 2\le q\le p+\nu$, if $$ \omega_{i}= L \lambda_{i}\|x_{i} - \hat{x}_{i-1}\|^{p+\nu-q}\le\theta_2\le 1,$$ Algorithm \ref{alg:dis} will converge according to Theorem \ref{thm:bound}; otherwise, the convergence behavior of Algorithm \ref{alg:dis} cannot be guaranteed. More specifically, when $q=p+\nu$ if $\omega_{i}$ satisfies \eqref{eq:lambda-cond}, then Algorithm \ref{alg:dis} converges according  to Theorem \ref{thm:bound-a}; when $2\le q < p+\nu,$ if $\omega_{i}$ satisfies \eqref{eq:lambda-cond-a}, then Algorithm \ref{alg:dis} converges according  to Theorem \ref{thm:bound-b}. When $q=p+\nu,$ we can trivially find $0<\theta_1\le\omega_{i}\le\theta_2\le 1$ to satisfy \eqref{eq:lambda-cond}; while when $2\le q <p+\nu,$ because $\omega_{i}$ involves the variable $x_{i}$ to optimize, it is nontrivial to find a $0<\theta_1\le\omega_{i}\le \theta_2<1$ to satisfy \eqref{eq:lambda-cond-a}.
A standard technique to ensure that the value of the convergence indicator $\omega_{i}$ stays in $[\theta_1, \theta_2]$ is through a binary search procedure \cite{monteiro2013accelerated,jiang2018optimal,bubeck2018near}. However, as per our discussion at the end of Section \ref{sec:dis}, the cost of the binary search procedure could substantially reduce the advantage of convergence rate of the optimal method in practice. 

\subsection{A Good Heuristic for Practical Implementation}
In this section, inspired by the analysis of Theorem \ref{thm:bound-b}, for the Algorithm \ref{alg:dis} with $2\le q<p+\nu$, instead of using a binary search,  we introduce a simple heuristic: in the $i$-th iteration of Algorithm \ref{alg:dis}, $A_{i}$ is set as its lower bound such that 
\begin{eqnarray}
A_{i} = \frac{C_0}{L}\big(h(x^*;x_0)\big)^{-\frac{p+\nu-q}{q}}\left(\frac{i}{p+\nu}\right)^{\frac{(q+1)(p+\nu)-q}{q}},\label{eq:heuristic}
\end{eqnarray}
where all the constants are from Theorem \ref{thm:bound-b}. With so assigned $A_i$,  $\lambda_{i}$ and $a_{i}$ can be easily determined by \eqref{eq:lam-A-x}. Therefore the per-iteration cost under the setting $2\le q<p+\nu$ will remain the same as the setting $q = p+\nu$.

However, if we use the heuristic  \eqref{eq:heuristic} of $A_{i}$ for $2\le q<p+\nu$, there is no theoretical guarantee for convergence of the algorithm. In this section, we conduct experiments to show that this heuristic  \eqref{eq:heuristic} is surprisingly effective: the values of the convergence indicator \eqref{eq:lambda-2} will always remain within the range $(0,1)$, hence Algorithm \ref{alg:dis} converges according to Theorem \ref{thm:bound-b}.

To be more precise, we consider the commonly second-order ($i.e., p=2$) setting with Euclidean Lipschitz smoothness Hessians ($i.e,$ $\nu =1$), and set $h(x;x_0):=\frac{1}{q}\|x-x_0\|_2^q$, where $q$ is chosen as $q\in\{2, 2.5, 3\}\subset [2, p+\nu]$. 
Meanwhile, as shown in Theorems \ref{thm:bound-a} and \ref{thm:bound-b}, the parameter $\alpha$ of Algorithm \ref{alg:dis} has only a minor influence on performance. Therefore to simplify our implementation, we always set $\alpha=1$. 
By setting $\alpha=1$, when $p=2, \nu=1,$  given $\hat{x}_{i-1}$ and $\lambda_{i}$, the subproblem of finding $x_{i}$  in the Step 5 of UAF is a standard cubic regularized Newton step \cite{carmon2016gradient}. We solve this subproblem to high accuracy by an implementaion \cite{kohler2017sub} \footnote{The GitHub URL: https://github.com/dalab/subsampled\_cubic\_regularization} of the Lanzcos method \cite{carmon2018analysis}.  Furthermore, in the heuristic \eqref{eq:heuristic} for $A_{i}$, $C_0$ is determined by the parameters $p, \nu, q, \theta_1, \theta_2, \beta$ and $\gamma$, while we already set the values of $p, \nu$. By the uniformly convexity of $h(x;x_0)=\frac{1}{q}\|x-x_0\|^q_2 (q\ge 2)$ \cite{nesterov2008accelerating}, we have $\gamma = \beta = 2^{2-q}$. We simply choose $\theta_1 = 0.5, \theta_2 = 0.67$. The Lipschitz smoothness constant $L$ is tuned in $\{10^{-3}, 10^{-2}, 10^{-1}, 1, 10,10^2, 10^3\}$ to optimize the convergence speed in terms of run time, while the value of $h(x^*;x_0)=\frac{1}{q}\|x^*-x_0\|_2^q$ is determined by setting $x_0 = 0$ and using an approximation of $x^*$ to replace $x^*$.

Under the above setting, three instances of the UAF Algorithm \ref{alg:dis} with $q=2, 2.5, 3$ respectively will be tested. The instance with $q=3$ is equivalent to the accelerated cubic regularized Newton method (ACNM) \cite{nesterov2008accelerating}. For the instance with $q=2$ or $2.5$, we always use the heuristic \eqref{eq:heuristic} to determine the values of $A_{i},$ $a_{i}$ and  $\lambda_{i}$ in each iteration.

\subsection{Experiments on Large-Scale Classification Datasets}
To verify the performance of the proposed UAF and the effectiveness of the heuristic \eqref{eq:heuristic} in all three instances, we consider large-scale optimization associated with the logistic regression problem as follows
\begin{eqnarray}
\min_{x\in\bbR^d}\left\{f(x):= \frac{1}{n} \sum_{j=1}^n \log(1+\exp(-\bar{b}_j \bar{a}_j^T x)) \right\},\label{eq:obj}
\end{eqnarray}
where  $\{(\bar{a}_j, \bar{b}_j)\}_{j=1}^n$ denotes a dataset. (For $j\in[n], \bar{a}_j\in\bbR^d$ denotes the $j$-th sample and $\bar{b}_j\in\{1,-1\}$ denotes the corresponding label of $\bar{a}_j$.) In our experiments, we choose the three datasets ``gisette\_scale'', ``a9a'' and ``w8a'' from the LIBSVM library \cite{CC01a} to validate the performance of our algorithm. 

\begin{figure}[t]
\centering
\subfigure{
\includegraphics[scale=0.21]{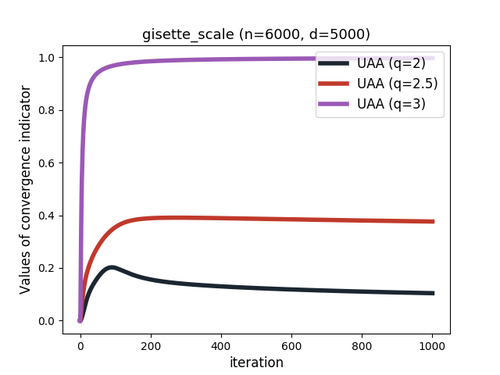}
}
\quad
\subfigure{
\includegraphics[scale=0.21]{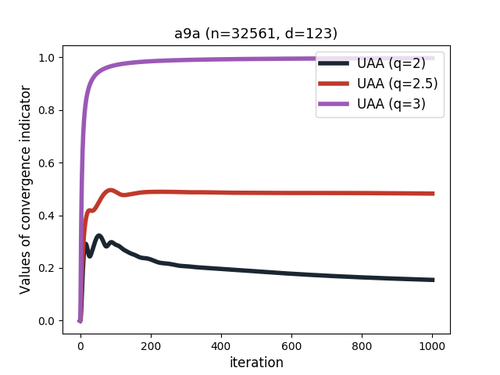}
}
\quad
\subfigure{
\includegraphics[scale=0.21]{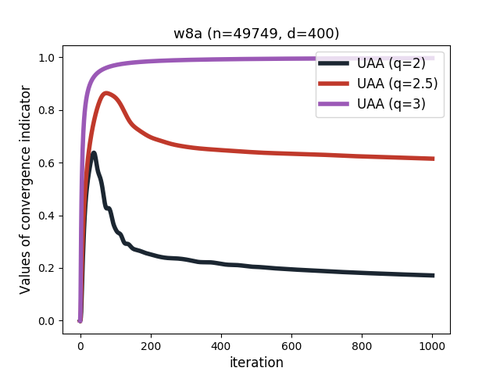}
}
\caption{The values of the convergence indicator \eqref{eq:lambda-2} versus the number of iterations, for  the $3$ datasets ``gisette\_scale'',  ``a9a'', ``w8a,'' respectively.}\label{fig:indicator}
\end{figure}

\begin{figure}[t]
\centering
\subfigure{
\includegraphics[scale=0.21]{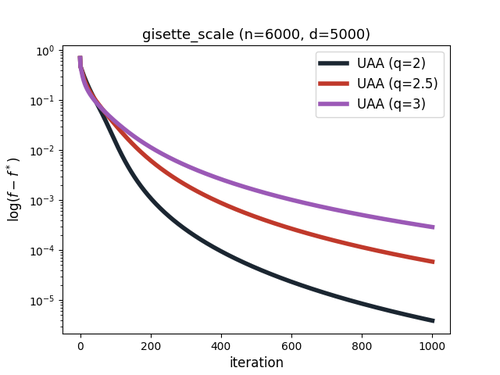}
}
\quad
\subfigure{
\includegraphics[scale=0.21]{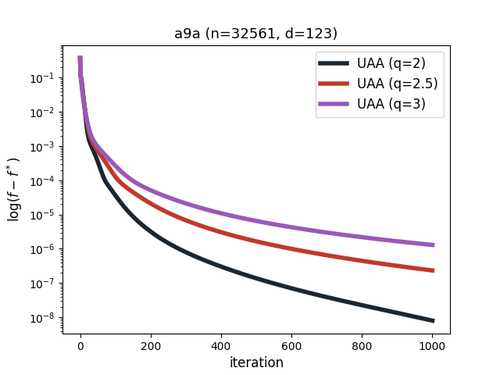}
}
\quad
\subfigure{
\includegraphics[scale=0.21]{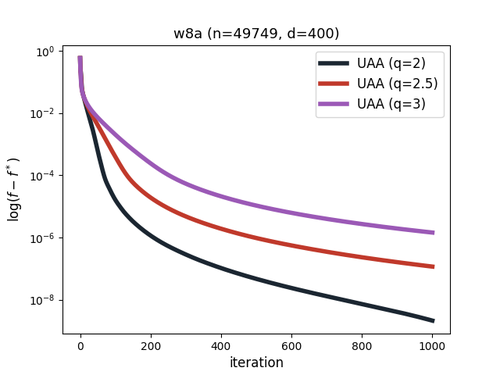}
}
\caption{Accuracy of  the objective function \eqref{eq:obj} versus the number of iterations, for the $3$ datasets ``gisette\_scale'',  ``a9a'', ``w8a,'' respectively.}\label{fig:iteration}
\end{figure}

\begin{figure}[t]
\centering
\subfigure{
\includegraphics[scale=0.21]{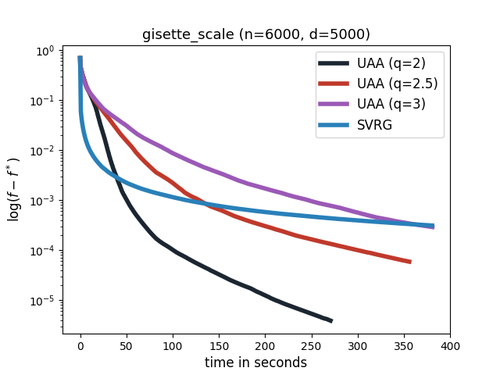}
}
\quad
\subfigure{
\includegraphics[scale=0.21]{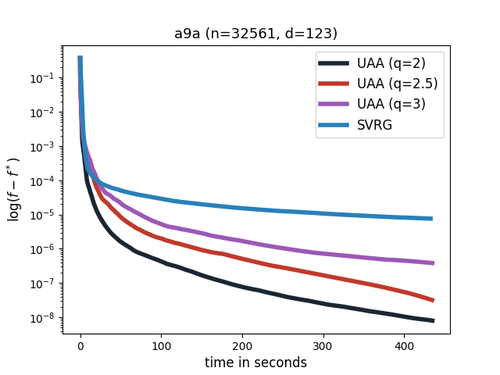}
}
\quad
\subfigure{
\includegraphics[scale=0.21]{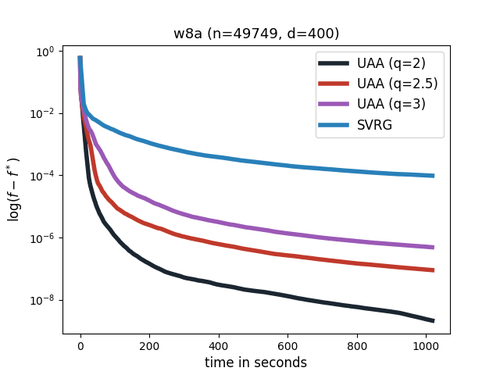}
}
\caption{Accuracy of the objective function \eqref{eq:obj} versus algorithm run time, for the $3$ datasets ``gisette\_scale'',  ``a9a'', ``w8a,'' respectively.}\label{fig:time}
\end{figure}

In Figure \ref{fig:indicator}, 
we show the values of the convergence indicator \eqref{eq:lambda-2} of UAF  along the iterations. It is interesting (and somewhat surprising) to see that after several initial steps, the convergence indicator will approach to a constant in $[0, 1]$. For the case with $q=3$, $i.e.,$ the ACNM \cite{nesterov2008accelerating}, the value of the indicator will approach to $1$, which matches the condition \eqref{eq:lambda-cond} with the optimal choice $\theta_1=\theta_2=1$.  For the cases with $q=2$ and $2.5$, the values of the indicator will stay stable around a constant in $(0,1).$ 

Because the values of the indicators satisfy the condition \eqref{eq:lambda-cond} when $q=3$ and the condition \eqref{eq:lambda-cond-a} when $q=2$ and $2.5$, the UAF algorithm will converge according to the rates in Theorems \ref{thm:bound-a} and \ref{thm:bound-b} respectively, which is shown in Figure \ref{fig:iteration}. In Figure \ref{fig:iteration}, with the heuristic \eqref{eq:heuristic}, then the UAF with $q=2$ has the fastest convergence speed, which matches the theoretical result that the setting $q=2$ gives us the best possible iteration complexity $O\left(k^{-\frac{3(p+\nu)-2}{2}}\right)$.

An interesting phenomenon is that the speed edge for the cases $q=2$ and  $2.5$ is beyond our expectation based on the bound \eqref{eq:bound-b}. In the $k$-th iteration, from the theoretical bound in Theorems \ref{thm:bound-a} and \ref{thm:bound-b}, the error ratio from the setting $q=3$ to the setting $q\in[2, p+\nu)$ should be
\begin{equation}
O\left(\frac{C_0}{\theta_1 c_q \gamma} \left(\frac{k}{(p+\nu)h(x^*;x_0)}\right)^{\frac{p+\nu-q}{q}} \right). \label{eq:ratio}
\end{equation}
In the experiments on all the $3$ datasets, we found empirically that $h(x^*;x_0)>1$. Meanwhile, by simple calculation, we also know that $\frac{\theta_1 c_q \gamma}{C_0}>1$. Therefore in the $1000$-th iteration, by the theoretical bound \eqref{eq:ratio}, the error ratio from $q=3$ to $2$ should not go beyond $(\frac{1000}{3})^{\frac{1}{2}}<{20}$. However, in practice the ratio is beyond ${100}$. A possible explanation for this phenomenon is that even we do not add any strongly convex regularizer in \eqref{eq:obj}, the problem itself may have some kind of local strong convexity around the minimum point (also known as implicit regularization). Such implicit strong convexity makes the algorithms converges faster as the iterate approaches the minimizer.

In Figure \ref{fig:time}, we show the performance comparison measured by error versus run time. Here we add a \emph{stochastic variance reduction gradient} (SVRG) \cite{johnson2013accelerating} method to show the practical efficiency of the proposed UAF algorithm.  SVRG is a representative first-order algorithm for finite-sum stochastic convex optimization. The implementation of SVRG is also from the GitHub project of \cite{kohler2017sub} and the learning rate of SVRG is tuned in $\{10^{-4}, 10^{-3}, 10^{-2}, 10^{-1}, $ $ 10, 10^2\}$. 

As shown in Figure \ref{fig:time}, SVRG can effectively exploit the finite-sum structure of the objective \eqref{eq:obj} and shows advantage in obtaining a low-accurate solution quickly. However, when further pursuing a high-accuracy solution, the high-order UAFs demonstrate clear edges of their faster convergence rates. In particular, with the effective heuristic \eqref{eq:heuristic}, the UAF with $q=2$ demonstrates consistent and superior performance in terms of run time behaviors.

\section{Conclusions}
In this paper, inspired by recent work on high-order acceleration methods, we have introduced a rather unified framework towards developing and understanding high-order acceleration algorithms for convex optimization. We show how various ideas, techniques, results, and algorithms can be derived from a simple vanilla proximal method (VPM). This perspective also helps reveal connections and clarify (rather significant) differences between the continuous setting and the discrete-time setting for acceleration. Based on this framework, through careful analysis, we are able to derive a unified acceleration framework (UAF) that achieves the optimal lower bounds for functions that have H\"{o}lder continuous derivatives. For functions that are uniformly convex, we have introduced a general restart scheme that helps our algorithm to achieve the optimal bounds. Our analysis and results also seem to unify many results known for the first order and high order methods, as well as results previously obtained through two separate approaches, namely the ACNM \cite{nesterov2008accelerating} and A-NPE \cite{monteiro2013accelerated} approaches. Meanwhile, the UAF is the first high-order acceleration approach that can be used in general (non-Euclidean) norm settings. Furthermore, for practical implementation of the proposed algorithm, through a new heuristic inspired from our analysis, our experiments show how the binary search procedure required by the optimal acceleration methods can be significantly simplified or forgone. This helps alleviate concerns about practical efficiency of optimal high-order acceleration methods versus suboptimal ones \cite{nesterov2018implementable}.

\appendix

\section{Proofs for Section \ref{sec:prelim}}\label{sec:append-prelim}

\subsection{Proof of Example \ref{exam:smooth}}\label{sec:exam:smooth}

$ $

\begin{proof}[Proof of Example \ref{exam:smooth}]
By direct computation, for $x, h\in\bbR^d$, we have 
\begin{eqnarray}
\langle \nabla^2 f(x) h, h\rangle& =& \frac{1}{n}\sum_{j=1}^n \frac{\exp(-\bar{b}_j \bar{a}_j^T x)}{ (1+ \exp(-\bar{b}_j \bar{a}_j^T x))^2}\langle \bar{a}_j, h\rangle^2\nonumber\\
&\le& \frac{1}{n}\sum_{j=1}^n\langle \bar{a}_j, h\rangle^2 = h^T B h.
\end{eqnarray}
Meanwhile, we have
\begin{eqnarray}
\nabla^3 f(x)[h, h, h]& = &\frac{1}{n}\sum_{j=1}^n -  \frac{\exp(-\bar{b}_j \bar{a}_j^T x)(1-\exp(-\bar{b}_j \bar{a}_j^T x))}{(1+\exp(-\bar{b}_j \bar{a}_j^T x))^3  }  \langle \bar{b}_j \bar{a}_j, h\rangle ^3  \; \le \; \frac{1}{n}\sum_{j=1}^n  |\langle \bar{a}_j, h\rangle|^3\nonumber\\
&\le& \frac{1}{n}(\sum_{j=1}^n  \langle \bar{a}_j, h\rangle ^2)\max_{j\in[n]} |\langle \bar{a}_j, h\rangle|\nonumber\\
&=& h^T B h\cdot \max_{j\in[n]} |\langle \bar{a}_j, h\rangle|.
\end{eqnarray}

Therefore,
\begin{eqnarray}
\|\nabla^2 f(x)\|_{q} = \max_{h\in\bbR^d:\|h\|_p\le 1} \langle \nabla^2 f(x)h,h\rangle\le\max_{h\in\bbR^d:\|h\|_p\le 1}\|Bh\|_q\le \|B\|_{p, q}.
\end{eqnarray}
\begin{eqnarray}
\|\nabla^3 f(x)\|_{q} = \max_{h\in\bbR^d:\|h\|_p\le 1} \nabla^3 f(x)[h, h, h]\le \|B\|_{p, q} \max_{j\in[n]}\|\bar{a}_j\|_q.
\end{eqnarray}
Let 
\begin{equation}
L(\nu):=\sup_{x, y\in\bbR^d, x\neq y}\frac{\|\nabla^2 f(x) - \nabla^2 f(y)\|_q}{\|x-y\|_p^{\nu}}, \nu\in[0,1].
\end{equation}
Then we have $L(0) = \|\nabla^2 f(x)\|_{q} , L(1)=\|\nabla^3 f(x)\|_{q}.$
Note that $L(\nu)$ is log-convex, therefore we have 
\begin{equation}
L(\nu)\le L^{1-\nu}(0)L^{\nu}(1)\le \|B\|_{p,q}\max_{j\in [n]}\|\bar{a}_j\|_q^{\nu}.
\end{equation}

Example \ref{exam:smooth} is proved.

\end{proof}

\subsection{Proof of Lemma \ref{lem:hat-tilde-f} }\label{sub:lem:hat-tilde-f}

$ $

$ $

\begin{proof}
[Proof of Lemma  \ref{lem:hat-tilde-f}]
By the convexity of $g(x)$, \eqref{eq:hat-prop} holds trivially.

If $g(x)$ has $p$-th derivatives, for $i\in\{0,1,2,\ldots, p-1\}$, we define a sequence
\begin{equation}
C_i:=\frac{1}{i!}\int_{0}^1(1-\tau)^{i}\nabla^{i+1} g(y+\tau(x-y))[x-y]^{i+1}\d\tau.  \label{eq:C-i-def} 
\end{equation}
Then one has 
\begin{eqnarray}
C_0 &=& \int_0^1  \nabla g(y+\tau(x-y))[x-y]\d\tau = \int_0^1 \langle \nabla g(y+\tau(x-y)), x-y\rangle\d\tau\nonumber\\
&=&  g(y+\tau(x-y))|_{\tau=0}^1\nonumber\\
&=& g(x) - g(y).\label{eq:hat-tilde-C-0}
\end{eqnarray}
Meanwhile,
\begin{eqnarray}
\quad C_i &= &\frac{1}{i!}\int_{0}^1(1-\tau)^{i}\d\left(\nabla^{i}g(y+\tau(x-y))[x-y]^i\right)\nonumber\\
&=&\frac{1}{i!}\left(\nabla^{i}g(y+\tau(x-y))[x-y]^i\right)(1-\tau)^{i}|_{\tau=0}^1 \nonumber\\ &&\quad-\frac{1}{i!}\int_{0}^1\left(\nabla^{i}g(y+\tau(x-y))[x-y]^i\right)\d (1-\tau)^{i}\nonumber\\
&=&-\frac{1}{i!}\nabla^{i}g(y)[x-y]^i + \frac{1}{(i-1)!}\int_{0}^1 (1-\tau)^{i-1} \left(\nabla^{i}g(y+\tau(x-y))[x-y]^i\right)\d \tau\nonumber\\
&=& -\frac{1}{i!}\nabla^{i}g(y)[x-y]^i + C_{i-1}.\label{eq:hat-tilde-C-i}
\end{eqnarray}
Therefore by \eqref{eq:hat-tilde-C-0} and \eqref{eq:hat-tilde-C-i},   one has
\begin{eqnarray}
C_{p-1} &=& \sum_{i=1}^{p-1}(C_i - C_{i-1}) + C_0=\sum_{i=1}^{p-1}-\frac{1}{i!}\nabla^i g(y)[x-y]^i+ g(x) - g(y)\nonumber\\
&\overset{\circlenum{1}}{=}&f(x) - \tilde{f}(x;y) + \frac{1}{p!}\nabla^p g(y)[x-y]^p\nonumber\\
&\overset{\circlenum{2}}{=}&f(x) - \tilde{f}(x;y) + \frac{1}{(p-1)!}\nabla^p g(y)[x-y]^p \int_0^1 (1-\tau)^{p-1}\d\tau,\label{eq:C-p-1}
\end{eqnarray}
where \circlenum{1} is by the definition of $f(x)$ in \eqref{eq:prob} and the definition of $\tilde{f}(x;y)$ in \eqref{eq:tilde-f}, \circlenum{2} is by the fact that $\int_{0}^1 (1-\tau)^{p-1}\d\tau = \frac{1}{p}$.

Then by \eqref{eq:C-p-1}, it follows that
\begin{eqnarray*}
&&|f(x) - \tilde{f}(x;y)| \\
&=& \left|C_{p-1}  - \frac{1}{(p-1)!}\nabla^p g(y)[x-y]^p \int_0^1 (1-\tau)^{p-1}\d\tau\right|\\
&\overset{\circlenum{1}}{=}&\frac{1}{(p-1)!}\left| \int_0^1 (1-\tau)^{p-1}\left( \nabla^{p} g(y+\tau(x-y)) -  \nabla^p g(y) \right)[x-y]^p\d\tau\right|\\
&\le& \frac{1}{(p-1)!} \int_0^1 (1-\tau)^{p-1}\d\tau\max_{\tau\in[0,1]}\left|\left( \nabla^{p} g(y+\tau(x-y)) -  \nabla^p g(y) \right)[x-y]^p\right|\\
&\overset{\circlenum{2}}{\le}& \frac{1}{(p-1)!} \int_0^1 (1-\tau)^{p-1}\d\tau\max_{\tau\in[0,1]}\left\|\nabla^{p} g(y+\tau(x-y)) -  \nabla^p g(y)\right\|_* \|x-y\|^p\\
&\overset{\circlenum{3}}{\le}& \frac{1}{(p-1)!} \frac{1}{p}\max_{\tau\in[0,1]}\left((p-1)!L\|\tau(x-y)\|^{\nu}\right)\|x-y\|^p\\
&\le&\frac{L}{p}\|x-y\|^{p+\nu},
\end{eqnarray*}
where \circlenum{1} is by \eqref{eq:C-i-def}, and \circlenum{2} is the definition of the operator norm $w.r.t.$ a norm $\|\cdot\|$ of symmetric $p$-linear form in \eqref{eq:operator-norm}, and  \circlenum{3} is by Definition \ref{def:holder-deri} and the fact $\int_0^1 (1-\tau)^{p-1}\d\tau = \frac{1}{p}$.
Therefore \eqref{eq:tilde-prop-1} holds.

By \eqref{eq:C-p-1}, by taking gradient $w.r.t.$ $x$,  one has
\begin{eqnarray}
&&\nabla C_{p-1} \nonumber\\
&=& \nabla f(x) - \nabla \tilde{f}(x;y)+\frac{1}{(p-1)!} \nabla^pg(y)[x-y]^{p-1}\nonumber\\
 &=&  \nabla f(x) - \nabla \tilde{f}(x;y)+\frac{p}{(p-1)!} \nabla^pg(y)[x-y]^{p-1}\int_0^1 (1-\tau)^{p-1}\d\tau,\label{eq:c-p-1-2}
\end{eqnarray}
while by \eqref{eq:C-i-def}, one also has
\begin{eqnarray}
\nabla C_{p-1} = \frac{p}{(p-1)!}\int_0^1(1-\tau)^{p-1}\nabla ^p  g(y+\tau(x-y))[x-y]^{p-1}\d\tau.\label{eq:c-p-1-3}
\end{eqnarray}
By \eqref{eq:c-p-1-2} and \eqref{eq:c-p-1-3}, it follows that
\begin{eqnarray*}
&&\|\nabla f(x) - \nabla \tilde{f}(x;y)\|_* = \left\|\nabla C_{p-1} - \frac{p}{(p-1)!} \nabla^pg(y)[x-y]^{p-1}\int_0^1 (1-\tau)^{p-1}\d\tau\right\|_*\\
&=&\left\| \frac{p}{(p-1)!}\int_0^1 (\nabla^pg(y+\tau(x-y)) - \nabla^pg(y))[x-y]^{p-1} (1-\tau)^{p-1}\d\tau\right\|_*\\
&\overset{\circlenum{1}}{=}&\max_{v:\|v\|\le 1}\frac{p}{(p-1)!}\int_0^1 (\nabla^pg(y+\tau(x-y)) - \nabla^pg(y))[v][x-y]^{p-1} (1-\tau)^{p-1}\d\tau    \\
&\le&\frac{p}{(p-1)!}\int_0^1 \max_{v:\|v\|\le 1} (\nabla^pg(y+\tau(x-y)) - \nabla^pg(y))[v][x-y]^{p-1} (1-\tau)^{p-1}\d\tau    \\
&\le&\frac{p}{(p-1)!}\int_0^1(1-\tau)^{p-1}\d\tau \max_{\tau\in[0,1]} \max_{v:\|v\|\le 1} (\nabla^pg(y+\tau(x-y)) - \nabla^pg(y))[v][x-y]^{p-1}     \\
&\overset{\circlenum{2}}{\le}& \frac{p}{(p-1)!}\cdot \frac{1}{p}\cdot\max_{\tau\in[0,1]}\max_{v:\|v\|\le 1} \|\nabla^pg(y+\tau(x-y)) - \nabla^pg(y)\|_*\cdot \|v\| \cdot \|x-y\|^{p-1}\\
&=& \frac{p}{(p-1)!}\cdot \frac{1}{p}\cdot\max_{\tau\in[0,1]}\|\nabla^pg(y+\tau(x-y)) - \nabla^pg(y)\|_*\cdot\|x-y\|^{p-1}\\
&\overset{\circlenum{3}}{\le}& \frac{p}{(p-1)!}\cdot \frac{1}{p}\cdot\max_{\tau\in[0,1]} (p-1)!L\|\tau(x-y)\|^{\nu} \cdot\|x-y\|^{p-1}\\
&\le&L\|x-y\|^{p+\nu-1},
\end{eqnarray*}
where \circlenum{1} is by the definition of $\|\cdot\|_*$, \circlenum{2}  is by the fact $\int_0^1 (1-\tau)^{p-1}\d\tau = \frac{1}{p}$ and the definition of the operator norm $w.r.t.$ a norm $\|\cdot\|$ of symmetric $p$-linear form in \eqref{eq:operator-norm}, and  \circlenum{3} is by Definition \ref{def:holder-deri}. 
Therefore \eqref{eq:tilde-prop-2} holds.

Lemma \ref{lem:hat-tilde-f} is proved.

\end{proof}

\subsection{Proof of Lemma \ref{lem:seq}}\label{sub:lem:seq}
$ $

\medskip

\begin{proof}[Proof of Lemma \ref{lem:seq}]

For \eqref{eq:seq-1}, one has 
\begin{eqnarray}
b_k-b_{k-1} \ge C^{\frac{1}{\rho}}b_k^{\frac{\rho-1}{\rho}}.\nonumber
\end{eqnarray}
Then by $b_0=0$, 
\begin{eqnarray}
b_k = \sum_{i=1}^k(b_i-b_{i-1}) \ge C^{\frac{1}{\rho}} \sum_{i=1}^k b_i^{\frac{\rho-1}{\rho}}.\nonumber
\end{eqnarray}
Then in Lemma \ref{eq:bjl}, for $ i \ge 1$, by setting $B_i := b_i^{\frac{\rho-1}{\rho}}, \upsilon := \frac{\rho}{\rho-1}, c := C^{\frac{1}{\rho}}$, then one has
\begin{eqnarray}
b_k^{\frac{\rho-1}{\rho}} = B_k \ge \left(\frac{\upsilon-1}{\upsilon}c k\right)^{\frac{1}{\upsilon-1}} = \left( \frac{1}{\rho}C^{\frac{1}{\rho}}k\right)^{\rho-1} .\nonumber
\end{eqnarray}
Then after a simple rearrangement, we obtain \eqref{eq:seq-2}.

\medskip\medskip

For \eqref{eq:seq-3}, 
 by using the reverse H\"{o}lder inequality, $\|fg\|_1\ge \|f\|_{\frac{1}{t}}\|g\|_{-\frac{1}{t-1}}$ for $t\ge 1$ and invoking this with $t = \rho\delta +1$ and by $b_0=0$, then 
\begin{eqnarray}
&&\sum_{i=1}^k \left(\frac{b_i^{\rho-1}}{(b_i-b_{i-1})^{\rho}}\right)^\delta = \sum_{i=1}^k b_i^{(\rho-1)\delta} (b_i-b_{i-1})^{-\rho\delta}\nonumber\\ 
&\ge&\left(\sum_{i=1}^k b_i^{(\rho-1)\delta\cdot\frac{1}{t}}\right)^t \left(\sum_{i=1}^k (b_i - b_{i-1})^{-\rho\delta\cdot\frac{-1}{t-1}}\right)^{-(t-1)}\nonumber\\ 
&=&\left(\sum_{i=1}^k b_i^{\frac{(\rho-1)\delta}{\rho\delta +1}}\right)^{\rho\delta +1}\left(\sum_{i=1}^k (b_i - b_{i-1})\right)^{-\rho\delta}\nonumber\\ 
&=&\left(\sum_{i=1}^k b_i^{\frac{(\rho-1)\delta}{\rho\delta +1}}\right)^{\rho\delta +1} b_k^{-\rho\delta}.
\end{eqnarray} 
Then by \eqref{eq:seq-3}, we have
\begin{eqnarray}
b_k^{\frac{ \rho\delta}{ \rho\delta+1} }\ge C^{-\frac{1}{\rho\delta+1}}&\left(\sum_{i=1}^k b_i^{\frac{(\rho-1)\delta}{\rho\delta +1}}\right).
\end{eqnarray}
Then in Lemma \ref{eq:bjl}, for $i\ge 1$, by setting $B_i := b_i^{\frac{(\rho-1)\delta}{\rho\delta+1}}, \upsilon := \frac{\rho}{\rho-1}, c := C^{-\frac{1}{\rho\delta+1}}$, then one has
\begin{eqnarray}
b_k^{\frac{(\rho-1)\delta}{\rho\delta+1}} =  B_k \ge \left(\frac{\upsilon-1}{\upsilon}c k\right)^{\frac{1}{\upsilon-1}}=\left( \frac{1}{\rho} C^{-\frac{1}{\rho\delta+1}} k\right)^{\rho-1}.
\end{eqnarray}
Then after a simple rearrangement, we obtain \eqref{eq:seq-4}.

Lemma \ref{lem:seq} is proved.

\end{proof}

\section{Proofs for Section \ref{sec:cont}}
\label{sec:cont-proofs}

\subsection{Proof of Lemma \ref{lem:cont-upper}}\label{sub:lem:cont-upper}
$ $

\begin{proof}[Proof of Lemma \ref{lem:cont-upper}] 
By Lemma \ref{lem:hat-tilde-f}, we have $\hat{f}({x;x_{\tau}})\le f(x)$. Thus one has
\[
\psi_t^{\rm{cont}}(x)\le \int_{0}^t a_{\tau} f(x)\d\tau + h(x; x_0)=A_t f(x)+h(x;x_0).
\]
Then it follows that
\[
\min_{x\in \bbR^d} \psi_t^{\rm{cont}}(x)\le \min_{x\in \bbR^d}\{ A_t f(x)+h(x;x_0)  \} \le A_t f(x^*)+h(x^*;x_0). 
\]

Lemma \ref{lem:cont-upper} is proved.
\end{proof}

\vspace{3cm}

\subsection{Proof of Lemma \ref{lem:cont}}\label{sub:lem:cont}
$ $

\begin{proof}[Proof of Lemma \ref{lem:cont}]
It follows that
\begin{eqnarray}
\quad\quad&&\frac{\d}{\d t}\psi_t^{\rm{cont}}(z_t) \nonumber\\
&=& \frac{\d}{\d t}\left(\int_{0}^t a_{\tau}\hat{f}(z_t; x_{\tau})\d\tau +h(z_t;x_0)\right)\nonumber\\
&{=}&a_t \hat{f}(z_t; x_t) +\int_{0}^t \langle a_{\tau}\nabla \hat{f}(z_t;x_{\tau}), \dot{z}_t\rangle\d\tau + \langle \nabla h(z_t;x_0), \dot{z}_{t}\rangle\nonumber\\
&=&a_t\hat{f}(z_t; x_t)+\left\langle \int_{0}^t a_{\tau} \nabla \hat{f}(z_t;x_{\tau})\d\tau+\nabla h(z_t;x_0), \dot{z}_{t}\right\rangle\nonumber\\
&\ge&a_t(\hat{f}(x_t;x_t)+\langle\nabla \hat{f}(x_t;x_t), z_t-x_t\rangle) +\left\langle \int_{0}^t  a_{\tau} \nabla \hat{f}(z_t;x_{\tau})\d\tau+\nabla h(z_t;x_0), \dot{z}_{t}\right\rangle\nonumber\\
&\overset{\circlenum{1}}{=}&a_t(f(x_t)+\langle\nabla f(x_t), z_t-x_t\rangle)+\left\langle \int_{0}^t a_{\tau} \nabla \hat{f}(z_t;x_{\tau})\d\tau+\nabla h(z_t;x_0), \dot{z}_{t}\right\rangle\nonumber
\end{eqnarray}
where \circlenum{1} is by the definition of $\hat{f}(x;y)$ in \eqref{eq:hat-f}.

Then by the optimality condition, one has 
$$\left\langle \int_{0}^t a_{\tau} \nabla \hat{f}(z_t;x_{\tau})\d\tau+\nabla h(z_t;x_0), \dot{z}_{t}\right\rangle\ge 0,$$
so 
\begin{eqnarray}
\frac{\d}{\d t}\psi_t^{\rm{cont}}(z_t)\ge a_t f(x_t)+ a_t\langle\nabla f(x_t), z_t-x_t\rangle.\label{eq:cont-tmp0}
\end{eqnarray}
Furthermore, one has
\begin{align}
\frac{\d (A_t f(x_t))}{\d t} = a_t f(x_t) + A_t \langle \nabla f(x_t), \dot{x}_t\rangle. \label{eq:cont-tmp1}
\end{align}

By Combing \eqref{eq:cont-tmp0} and \eqref{eq:cont-tmp1}, one has
\begin{equation}
\frac{\d}{\d t}\left(A_t f(x_t) - \min_{x\in \bbR^d}\psi_t^{\rm{cont}}(x)\right)\!=\!\frac{\d}{\d t}(A_t f(x_t) - \psi_t^{\rm{cont}}(z_t)) \le \langle \nabla f(x_t), A_t \dot{x}_t - a_t(z_t-x_t)\rangle.\label{eq:deri}
\end{equation}
Finally, by $A_0=0$ and $\min_{x\in \bbR^d}\psi_0(x)=0$, one has
\begin{eqnarray}
A_0 f(x_0) - \min_{x\in \bbR^d}\psi_0^{\rm{cont}}(x)=0. \label{eq:init}
\end{eqnarray}
By combing \eqref{eq:deri} and \eqref{eq:init}, and taking integral from $\tau=0$ to $1$, then Lemma \ref{lem:cont} is proved.

\end{proof}

\section{Proofs for Section \ref{sec:dis}}

\subsection{Proof of Lemma \ref{lem:dis-upper}}
\label{subsec:lem:dis-upper}
$ $

\begin{proof}[Proof of Lemma \ref{lem:dis-upper}]
By Lemma \ref{lem:hat-tilde-f}, we have $\hat{f}(x;x_i)\le f(x)$. Then one has 
\begin{align*}
\psi_k^{\text{dis}}(x)&\le \sum_{i=1}^{k}a_{i} f(x) + h(x; x_0)\\ 
&= A_kf(x)+h(x;x_0).
\end{align*}
Then it follows that
\[
\min_{x\in \bbR^d} \psi_k^{\rm{dis}}(x)\le \min_{x\in \bbR^d}\{ A_k f(x)+h(x;x_0)  \} \le A_k f(x^*)+h(x^*;x_0). 
\] 

Lemma \ref{lem:dis-upper} is proved.
\end{proof}

\subsection{Proof of Lemma \ref{lem:dis-lower}}
\label{subsec:lem:dis-lower}
$ $

\begin{proof}[Proof of Lemma \ref{lem:dis-lower}]
First, in \eqref{eq:dis}, by $A_0 =0$ and $z_0=x_0,$ we have
\begin{equation}
A_0 f(x_0) - \psi_0^{\rm dis}(z_0)=0.\label{eq:A-0-x-0}
\end{equation}
Then by our assumption, $\hat{f}(x; x_i)$ is convex  $w.r.t.$ $\|\cdot\|$ and $h(x;x_0)$ is $(q, \gamma)$-uniformly convex $w.r.t.$ $\|\cdot\|$. Therefore for all $x, y\in \bbR^d$, it follows that 
\begin{align}
\psi^{\rm{dis}}_{i}(x)\ge \psi^{\rm{dis}}_{i}(y)+\langle \nabla \psi^{\rm{dis}}_{i}(y), x-y\rangle  + \frac{ \gamma }{q}\|x-y\|^q.
\end{align}
Then by the optimality condition of $z_i$, it follows that for all $x\in \bbR^d$, $\langle \nabla \psi^{\rm{dis}}_{i}(z_i), x-z_i\rangle \ge 0$. Therefore, 
it follows that
\begin{eqnarray}
\psi^{\rm{dis}}_{i}(x)\ge \psi^{\rm{dis}}_i(z_i)  + \frac{ \gamma }{q}\|x-z_i\|^q.\label{eq:dis-tmp1}
\end{eqnarray}
Therefore by \eqref{eq:dis-tmp1}, 
\begin{eqnarray}
\psi^{\rm{dis}}_{i}(x) &=& \psi^{\rm{dis}}_{i-1}(x)+ a_{i}\hat{f}(x; x_{i})\nonumber\\
&\ge& \psi^{\rm{dis}}_{i-1}(z_{i-1}) + \frac{ \gamma }{q}\|x-z_{i-1}\|^q + a_{i}\hat{f}(x;x_{i}).\label{eq:dis-tmp2}
\end{eqnarray}
Meanwhile, we can lower bound the last term of RHS of \eqref{eq:dis-tmp2}.
\begin{eqnarray*}
a_{i}\hat{f}(x;x_{i}) &\overset{\circlenum{1}}{\ge}& a_{i}(\hat{f}(x_{i};x_{i}) + \langle \nabla \hat{f}(x_{i};x_{i}), x- x_{i}\rangle) \\
&\overset{\circlenum{2}}{=}&  a_{i}(f(x_{i}) + \langle \nabla f(x_{i}), x-x_{i}\rangle) \\ 
&\overset{\circlenum{3}}{=}& A_{i}(f(x_{i}) + \langle \nabla f(x_{i}), x-x_{i}\rangle) - A_{i-1}(f(x_{i}) + \langle \nabla f(x_{i}), x-x_{i}\rangle) \\
&=&A_{i}\left(f(x_{i}) + \left\langle \nabla f(x_{i}), \frac{a_{i}}{A_{i}}x+\frac{A_{i-1}}{A_{i}}x_{i-1}-x_{i}\right\rangle\right) \\
&&\quad-A_{i-1}(f(x_{i}) + \langle \nabla f(x_{i}), x_{i-1}-x_{i}\rangle) \\
&\overset{\circlenum{4}}{\ge}&A_{i}\left(f(x_{i}) + \left\langle \nabla f(x_{i}), \frac{a_{i}}{A_{i}}x+\frac{A_{i-1}}{A_{i}}x_{i-1}-x_{i}\right\rangle\right)-A_{i-1} f(x_{i-1}) \\
&=&A_{i}f(x_{i})-A_{i-1} f(x_{i-1}) + A_{i}  \left\langle \nabla f(x_{i}), \frac{a_{i}}{A_{i}}x+\frac{A_{i-1}}{A_{i}}x_{i-1}-x_{i}\right\rangle,
\end{eqnarray*}
where \circlenum{1} is by the convexity of $\hat{f}(x;y)$ $w.r.t.$ $x$, \circlenum{2} is by the definition of $\hat{f}(x;y)$ in \eqref{eq:hat-f},  \circlenum{3} is by the identity $a_{i} = A_{i} - A_{i-1}$, and \circlenum{4} is by the convexity of $f(x)$.  
 
Therefore it follows that
\begin{eqnarray}
\psi^{\rm{dis}}_{i}(x) &\ge& \psi^{\rm{dis}}_{i-1}(z_{i-1}) +  \frac{ \gamma }{q}\|x-z_{i-1}\|^q + A_{i}f(x_{i})-A_{i-1} f(x_{i-1})\nonumber\\
 &&+ A_{i}  \left\langle \nabla f(x_{i}), \frac{a_{i}}{A_{i}}x+\frac{A_{i-1}}{A_{i}}x_{i-1}-x_{i}\right\rangle. \label{eq:es-2-1}
\end{eqnarray}
By setting $x:=z_{i}$ and a simple arrangement of \eqref{eq:es-2-1}, we have
\begin{eqnarray}
&(A_{i}f(x_{i})- \psi^{\rm{dis}}_{i}(z_{i}))-(A_{i-1} f(x_{i-1})- \psi^{\rm{dis}}_{i-1}(z_{i-1}))\nonumber\\
&\le A_{i}\left\langle \nabla f(x_{i}), x_{i}-\frac{a_{i}}{A_{i}}z_{i}-\frac{A_{i-1}}{A_{i}}x_{i-1}\right\rangle -  \frac{ \gamma }{q}\|z_{i}-z_{i-1}\|^q \label{eq:es-2-2}
\end{eqnarray}
Summing \eqref{eq:es-2-2} from $i=0$ to $k-1$ and by \eqref{eq:A-0-x-0}, it follows that
\begin{eqnarray}
&&A_{k}f(x_{k})- \psi^{\rm{dis}}_{k}(z_{k})\nonumber\\
&\le& A_{0}f(x_{0})- \psi^{\rm{dis}}_{0}(z_{0}) \nonumber\\
&&+ \sum_{i=1}^{k}\left(A_{i}  \left\langle \nabla f(x_{i}), x_{i}-\frac{a_{i}}{A_{i}}z_{i}-\frac{A_{i-1}}{A_{i}}x_{i-1}\right\rangle - \frac{ \gamma }{q}\|z_{i}-z_{i-1}\|^q\right)\nonumber\\
&=&\sum_{i=1}^{k}\left(A_{i}  \left\langle \nabla f(x_{i}), x_{i}-\frac{a_{i}}{A_{i}}z_{i}-\frac{A_{i-1}}{A_{i}}x_{i-1}\right\rangle -\frac{ \gamma }{q}\|z_{i}-z_{i-1}\|^q\right).\nonumber
\end{eqnarray}

Then by the definition of $E_i$, Lemma \ref{lem:dis-lower} is proved.

\end{proof}

\subsection{Proof of Lemma \ref{lem:E-2}}
\label{subsec:lem:E-2}
$ $

\begin{proof}[Proof of Lemma \ref{lem:E-2}]

\medskip
By the definition of $E_{i}$, one has
\begin{eqnarray}
E_{i} &\overset{\circlenum{1}}{\le}& a_{i}\left\langle \nabla f(x_{i}),\hat{z}_{i}-z_{i}\right\rangle -  \frac{ \gamma }{q}\|z_{i}-z_{i-1}\|^q\nonumber\\
&\overset{\circlenum{2}}{\le}& a_{i}\left\langle \nabla f(x_{i}),\hat{z}_{i}-z_{i}\right\rangle -  \frac{ \gamma^{\prime}_{i} }{q}\|z_{i}-z_{i-1}\|^q\nonumber\\
&\overset{\circlenum{3}}{\le}&\left\langle  a_{i}\nabla f(x_{i}) + \gamma^{\prime}_{i}\nabla \frac{1}{q} \|\hat{z}_{i} - z_{i-1}\|^q,\hat{z}_{i}-z_{i}\right\rangle   \nonumber\\
&&\quad- \gamma^{\prime}_{i} \left( \frac{1}{q}\|\hat{z}_{i} - z_{i-1}\|^q + \frac{\beta}{q}\|\hat{z}_{i} - {z}_{i}\|^q\right)\nonumber\\
&\overset{\circlenum{4}}{=}& \left\langle a_{i} \nabla f(x_{i}) + \frac{\gamma^{\prime}_{i} A_{i}^{q-1}}{a_{i}^{q-1}}\nabla \frac{1}{q}\|x_{i} - \hat{x}_{i-1}\|^q,\hat{z}_{i}-z_{i}\right\rangle   \nonumber\\
&&\quad- \gamma^{\prime}_{i} \left( \frac{A_{i}^{q}}{q a_{i}^{q}}\|x_{i} - \hat{x}_{i-1}\|^q + \frac{\beta}{q}\|\hat{z}_{i} - {z}_{i}\|^q\right)\nonumber
\end{eqnarray}
where \circlenum{1} is by the definition of $\hat{z}_{i}$ in \eqref{eq:two-fixed}, \circlenum{2} is by the assumption that $\gamma \ge \gamma^{\prime}_{i}$, \circlenum{3} is by Assumption \ref{ass:norm} of $\frac{1}{q}\|\cdot\|^q$, \circlenum{4} is by the definition of $\hat{x}_{i-1}$ in \eqref{eq:two-fixed} such that $\nabla\|\hat{z}_{i} - z_i\|^q = \frac{A_{i}^{q-1}}{a_{i}^{q-1}}\nabla \|x_{i}-\hat{x}_{i-1}\|^q$. 

Lemma \ref{lem:E-2} is proved.

\end{proof}

\medskip\medskip\medskip\medskip

\medskip\medskip\medskip\medskip

\subsection{Proof of Lemma \ref{lem:E-3}}\label{sec:E-3} 
$ $

\begin{proof}[Proof of Lemma \ref{lem:E-3}]
By Lemma \ref{lem:E-2}, one has
\begin{eqnarray}
E_{i} &\le&  a_{i}\left\langle\nabla f(x_{i}) + \frac{\gamma^{\prime}_{i} A_{i}^{q-1}}{a_{i}^q}\nabla\frac{1}{q} \|x_{i} -\hat{x}_{i-1}\|^q,  \hat{z}_{i} - z_{i}\right\rangle \nonumber\\
&&-\gamma^{\prime}_{i}\left(
\frac{A_{i}^q}{qa_{i}^q}\|x_{i} -\hat{x}_{i-1}\|^q + \frac{\beta}{q}\|\hat{z}_{i} - z_{i}\|^q
  \right)\nonumber \\
&\le&a_{i}\langle \nabla f(x_{i}) - \nabla \tilde{f}(x_{i};\hat{x}_{i-1}), \hat{z}_{i} - z_{i}\rangle\nonumber\\
 &&+ a_{i}\left\langle \nabla \tilde{f}(x_{i};\hat{x}_{i-1}) + \frac{\gamma^{\prime}_{i} A_{i}^{q-1}}{a_{i}^q}\nabla\frac{1}{q} \|x_{i} -\hat{x}_{i-1}\|^q,  \hat{z}_{i} - z_{i}\right\rangle \nonumber\\
 &&-\gamma^{\prime}_{i}\left(
\frac{A_{i}^q}{qa_{i}^q}\|x_{i} -\hat{x}_{i-1}\|^q + \frac{\beta}{q}\|\hat{z}_{i} - z_{i}\|^q
  \right).\label{eq:E-3-1}
\end{eqnarray}
Meanwhile, it follows that
\begin{eqnarray}
&&a_{i}\langle \nabla f(x_{i}) - \nabla \tilde{f}(x_{i};\hat{x}_{i-1}), \hat{z}_{i} - z_{i}\rangle\nonumber\\
&&-\gamma^{\prime}_{i}\left(
\frac{A_{i}^q}{qa_{i}^q}\|x_{i} -\hat{x}_{i-1}\|^q + \frac{\beta}{q}\|\hat{z}_{i} - z_{i}\|^q
  \right)\nonumber \\
&\le&a_{i}\| \nabla f(x_{i}) - \nabla \tilde{f}(x_{i};\hat{x}_{i-1})\|_*\|\hat{z}_{i} - z_{i}\|\nonumber\\
 &&-\gamma^{\prime}_{i}\left(
\frac{A_{i}^q}{qa_{i}^q}\|x_{i} -\hat{x}_{i-1}\|^q + \frac{\beta}{q}\|\hat{z}_{i} - z_{i}\|^q
  \right)\nonumber \\
&\overset{\circlenum{1}}{\le}&a_{i}L\|x_{i} -\hat{x}_{i-1}\|^{p+\nu-1}\|\hat{z}_{i} - z_{i}\|\nonumber\\
 &&-\gamma^{\prime}_{i}\left(
\frac{A_{i}^q}{qa_{i}^q}\|x_{i} -\hat{x}_{i-1}\|^q + \frac{\beta}{q}\|\hat{z}_{i} - z_{i}\|^q
  \right)\nonumber \\
&\overset{\circlenum{2}}{\le}&\frac{q-1}{q}(\beta\gamma^{\prime}_{i})^{-\frac{1}{q-1}}(a_{i}L)^{\frac{q}{q-1}}\|x_{i}-\hat{x}_{i-1}\|^{\frac{q(p+\nu-1)}{q-1}}-\frac{\gamma^{\prime}_{i} A_{i}^q}{q a_{i}^q}\|x_{i}-\hat{x}_{i-1}\|^q\nonumber\\
&\overset{\circlenum{3}}{=}&\left(\left(L\frac{a_{i}^q}{c_q \gamma^{\prime}_{i} A_{i}^{q-1}} \|x_{i}-\hat{x}_{i-1}\|^{p+\nu-q}\right)^{\frac{q}{q-1}}   -1  \right)\frac{\gamma^{\prime}_{i} A_{i}^q}{q a_{i}^q}\|x_{i} -\hat{x}_{i-1}\|^q\nonumber\\
&\overset{\circlenum{4}}{=}& \left(\left(L\lambda_{i}^{\prime} \|x_{i}-\hat{x}_{i-1}\|^{p+\nu-q}\right)^{\frac{q}{q-1}} -1  \right)\frac{\gamma^{\prime}_{i} A_{i}^q}{q a_{i}^q}\|x_{i} -\hat{x}_{i-1}\|^q,\label{eq:E-3-2}
\end{eqnarray}
where \circlenum{1} is by \eqref{eq:E-4}, \circlenum{2} is by the technical Lemma \ref{lem:s-t}, \circlenum{3} is by
a simple rearrangement and the definition of $c_q$ in Lemma \ref{lem:E-3}, and \circlenum{4} is by the definition of $\lambda_{i}^{\prime}$.

Combing \eqref{eq:E-3-1} and \eqref{eq:E-3-2}, Lemma \ref{lem:E-3} is proved.

\end{proof}

\medskip\medskip\medskip\medskip

\subsection{Proof of Theorem \ref{thm:bound}}\label{sec:bound}
$ $

\begin{proof}[Proof of Theorem \ref{thm:bound}]
\medskip
First, for $i\ge 1$, if the conditions \eqref{eq:E-7} and \eqref{eq:lambda-2} are true, then one can know that \eqref{eq:E-6-2} and \eqref{eq:E-6-1} are true and thus for $i\ge 1$, $E_{i}\le 0$. Then by Lemma 
\ref{lem:dis-lower}, one has
\begin{eqnarray}
 A_kf(x_k) - \psi^{\rm{dis}}_k(z_k)\le \sum_{i=1}^{k}E_{i} \le 0.
 \end{eqnarray} 
Then combing Lemma \ref{lem:dis-upper}, one has
\begin{eqnarray}
 A_kf(x_k)\le \psi_k^{\rm{dis}}(z_k)\le A_k f(x^*)+h(x^*;x_0).
\end{eqnarray}

Theorem \ref{thm:bound} is proved.
\end{proof}

\subsection{Proof of Theorem \ref{thm:bound-a}}\label{sub:thm:bound-a}
$ $

\begin{proof}[Proof of Theorem \ref{thm:bound-a}]
First, by our assumption,   $ \{\lambda_{i}\}$ defined in \eqref{eq:lambda-1}	satisfies \eqref{eq:lambda-cond},  therefore $ \{\lambda_{i}\}$ satisfies \eqref{eq:lambda-2}; meanwhile $\{x_{i}\}$ satisfies \eqref{eq:E-7}. Therefore Theorem \ref{thm:bound} holds, $i.e.,$
\begin{eqnarray}
f(x_k) - f(x^*)\le \frac{h(x^*;x_0)}{A_k}.\label{eq:thm-3-1}
\end{eqnarray}
Then by \eqref{eq:lambda-cond}, because $L\lambda_{i} = \frac{La_{i}^q}{c_q\gamma A^{q-1}_{i}} =  \frac{L(A_{i}-A_{i-1})^q}{c_q\gamma A^{q-1}_{i}} \ge \theta_1$, in Lemma \ref{lem:seq}, by setting $b_i := A_i, \rho := p+\nu, C:=\frac{\theta_1 c_q \gamma}{L}$,  we can obtain the lower bound 
\begin{equation}
A_k\ge \frac{\theta_1 c_q \gamma}{L} \left(\frac{k}{p+\nu}\right)^{p+\nu}. \label{eq:A-k-a}   
\end{equation}
By combing \eqref{eq:thm-3-1} and \eqref{eq:A-k-a}, \eqref{eq:bound-a} is obtained.

Theorem \ref{thm:bound-a} is proved.
\end{proof}

\medskip\medskip\medskip\medskip

\subsection{Proof of Lemma \ref{lem:q-neq-p-nu}}\label{sub:lem:q-neq-p-nu}
$ $

\begin{proof}[Proof of Lemma \ref{lem:q-neq-p-nu}]
When \eqref{eq:E-7} and \eqref{eq:lambda-2} are satisfied, by Lemma \ref{lem:E-3}, we have
\begin{eqnarray}
\quad\quad\; E_{i}&\overset{\circlenum{1}}{\le}&(\theta_2^{\frac{q}{q-1}} -1) \frac{\gamma^{\prime}_{i} A_{i}^q}{q a_{i}^q}\|x_{i} -\hat{x}_{i-1}\|^q\nonumber\\
&\overset{\circlenum{2}}{\le}&(\theta_2^{\frac{q}{q-1}} -1) \frac{A_{i}^q}{q a_{i}^q}\left(\frac{\omega_{i}}{\theta_2}\right)^{\alpha}\gamma      \|x_{i} -\hat{x}_{i-1}\|^q\nonumber\\
&\overset{\circlenum{3}}{=}& (\theta_2^{\frac{q}{q-1}} -1) \frac{A_{i}^q}{q a_{i}^q}\left(\frac{\omega_{i}}{\theta_2}\right)^{\alpha}\gamma   (L\lambda_{i}\|x_{i} -\hat{x}_{i-1}\|^{p+\nu-q})^{\frac{q}{p+\nu-q}}(L\lambda_{i})^{-\frac{q}{p+\nu-q}}\nonumber\\
&\overset{\circlenum{4}}{=}&\frac{1}{q\theta_2^{\alpha}}(\theta_2^{\frac{q}{q-1}} -1)\frac{A_{i}^q}{ a_{i}^q}\omega_{i}^{\frac{\varsigma}{p+\nu-q}}\gamma(L\lambda_{i})^{-\frac{q}{p+\nu-q}}\nonumber\\
&\overset{\circlenum{5}}{=}&\frac{1}{q\theta_2^{\alpha}} (\theta_2^{\frac{q}{q-1}} -1)\omega_{i}^{\frac{\varsigma}{p+\nu-q}}\gamma
\left(\frac{c_q \gamma A_{i}^{p+\nu-1}}{L (A_{i} -A_{i-1})^{p+\nu}}\right)^{\frac{q}{p+\nu-q}},\label{eq:q-neq-p-nu-1}
\end{eqnarray}
where \circlenum{1} is by \eqref{eq:E-7} and \eqref{eq:lambda-2}, \circlenum{2} is by the value of $\gamma_{i}^{\prime}$ in \eqref{eq:gamma-prime-value} and the definition of $\omega_{i}$ in Lemma \ref{lem:q-neq-p-nu}, \circlenum{3} is by  a simple rearrangement, \circlenum{4} is by definition of $\omega_{i}$ and $\varsigma=\alpha(p+\nu) + (1-\alpha)q$, \circlenum{5} is by definition of $\lambda_{i}$ in \eqref{eq:lambda-1} and the fact $a_{i} = A_{i} - A_{i-1}.$

Then by combing Lemmas \ref{lem:dis-upper} and \ref{lem:dis-lower}, it follows that 
\begin{eqnarray}
A_k f(x_k) \le \psi_k^{\rm{dis}}(z_k) + \sum_{i=1}^{k} E_{i} \le  A_k f(x^*) + h(x^*;x_0) + \sum_{i=1}^{k} E_{i}.\label{eq:q-neq-p-nu-2}
\end{eqnarray}
Then by combing \eqref{eq:q-neq-p-nu-1} and \eqref{eq:q-neq-p-nu-2}, and $f(x_k)\ge f(x^*), A_k\ge 0$, one has
\begin{eqnarray*}
\frac{1}{q\theta_2^{\alpha}} (1-\theta_2^{\frac{q}{q-1}})\gamma\sum_{i=1}^{k}
\omega_{i}^{\frac{\varsigma}{p+\nu-q}}\left(\frac{c_q \gamma A_{i}^{p+\nu-1}}{L (A_{i} -A_{i-1})^{p+\nu}}\right)^{\frac{q}{p+\nu-q}}\le \sum_{i=1}^{k}-E_{i}\le h(x^*;x_0).
\end{eqnarray*}
Then after a simple rearrangement, Lemma  \ref{lem:q-neq-p-nu} is proved.
\end{proof}

\medskip\medskip\medskip\medskip
\subsection{Proof of Theorem \ref{thm:bound-b}}\label{sub:thm:bound-b}
$ $

\begin{proof}[Proof of Theorem \ref{thm:bound-b}]
First, by our assumption, $\{\lambda_{i}\}$ defined in \eqref{eq:lambda-1}	satisfies \eqref{eq:lambda-cond-a},  therefore $ \{\lambda_{i}\}$ satisfies \eqref{eq:lambda-2}; meanwhile $\{x_{i}\}$ satisfies \eqref{eq:E-7}. Therefore Theorem \ref{thm:bound} holds, $i.e.,$
\begin{eqnarray}
f(x_k) - f(x^*)\le \frac{h(x^*;x_0)}{A_k}.\label{eq:A-k-b}
\end{eqnarray}

Then by Lemma \ref{lem:q-neq-p-nu} and the assumption that $\omega_{i}\ge \theta_1$, we have
\begin{eqnarray}
\sum_{i=1}^{k}
\left(\frac{A_{i}^{p+\nu-1}}{(A_{i} -A_{i-1})^{p+\nu}}\right)^{\frac{q}{p+\nu-q}}\le   (C_0^{-1}L)^{\frac{q}{p+\nu-q}}h(x^*;x_0), 
\end{eqnarray}
where $C_0$ is defined in Theorem \ref{thm:bound-b}. 

In Lemma \ref{lem:seq}, for $1\le i\le k$, by setting $b_i:=A_i, \rho := p+\nu, \delta:= \frac{q}{p+\nu-q}, C:=(C_0^{-1}L)^{\frac{q}{p+\nu-q}}h(x^*;x_0)$, then we obtain the lower bound
\begin{eqnarray}
A_k&\ge \frac{C_0}{L}
\left(h(x^*;x_0)\right)^{-\frac{p+\nu-q}{q}}\left(\frac{k}{p+\nu}\right)^{ \frac{(q+1)(p+\nu)-q}{q}}.
\end{eqnarray}
Then combing \eqref{eq:A-k-b}, we obtain \eqref{eq:bound-b}.
\end{proof}

\medskip\medskip\medskip\medskip

\subsection{Proof of Proposition \ref{prop:exist}}\label{sec:prop:exist}
$ $

\begin{proof}[Proof of Proposition \ref{prop:exist}]
First by our assumption about $\|\cdot\|$ and $\tilde{f}(x;y)$, \eqref{eq:y-z} is a strictly convex function, therefore $w(\upsilon)$ is a continuous function of $\upsilon$. Meanwhile $x(\lambda)$ is continuous about $\lambda$. Therefore $\chi(\lambda)$ is continuous $w.r.t.$ $\lambda$.

Next by the fact $\varsigma=\alpha(p+\nu)+(1-\alpha)q\in[q, p+\nu]$ and 
\begin{eqnarray}
&&\tilde{f}(z;\upsilon) +   \frac{L^{\alpha}}{ q c_q  \lambda^{(1-\alpha)}\theta_2^{\alpha}\varsigma }   \|z-\upsilon\|^{\varsigma}
\le\tilde{f}(\upsilon;\upsilon)=f(\upsilon)<+\infty,
\end{eqnarray}
as $\lambda\rightarrow 0$, $\|z-\upsilon\|\rightarrow 0$ if $\varsigma\in[q,p+\nu)$ or $\|z-\upsilon\|$ is a finite value if $\varsigma=p+\nu$. In both cases, we have $\chi(0) = 0$. Then since $f(\upsilon)\neq f(x^*)$, we will also have
as $\lambda\rightarrow+\infty$, it is easy to find that $\frac{a(\lambda)}{A+a(\lambda)}\rightarrow 1$ and thus $x(\lambda) = x$. Since $f(x)\neq f(x^*)$, we have $\omega(x)\neq x$.  Therefore $\chi(+\infty) = +\infty.$ 
\end{proof}

\medskip\medskip\medskip

\medskip\medskip\medskip\medskip

\medskip

In Lemma \ref{lem:s-t}, we give a simplified version of the inequality in \cite[Lemma 2]{nesterov2008accelerating}.
\begin{lemma}[\cite{nesterov2008accelerating}]\label{lem:s-t}
For $s, t\in \bbR_+$ and $q\ge 2, \sigma>0$, one has
\begin{eqnarray}
|st|\le \frac{\sigma}{q}t^q + \frac{q-1}{q}\left(\frac{1}{\sigma}\right)^{\frac{1}{q-1}}s^{\frac{q}{q-1}}.
\end{eqnarray}
\end{lemma}

In Lemma \ref{eq:bjl}, we give  Lemma \ref{eq:bjl} from \cite{bubeck2018near}, which is used to prove Lemma \ref{lem:seq}.
\begin{lemma}[\cite{bubeck2018near}]\label{eq:bjl}
Given a positive sequence $\{B_i\}$ such that $B_k^{\upsilon} \ge c\cdot\sum_{i=1}^k B_i$, where $c>0, \upsilon >1$. Then it follows that
\begin{eqnarray*}
B_k \ge \left(\frac{\upsilon-1}{\upsilon}ck\right)^{\frac{1}{\upsilon-1}}.
\end{eqnarray*}
\end{lemma}

\medskip\medskip

\section{Proofs for Section \ref{sec:restart}}

\subsection{Proof of Theorem \ref{thm:restart}}\label{sec:thm:restart}
$ $

\begin{proof}[Proof of Theorem \ref{thm:restart}]

\medskip\medskip
\textbf{The case with $k\le k_0$.}

\medskip

Denote $R_k := R\cdot(\frac{1}{2})^k$. For $k_0-1 \ge k\ge 0,$ let us prove by induction that 
\begin{eqnarray}
\|y_k - x^*\|\le R_k. \label{eq:y-k-R}
\end{eqnarray}

By the optimality of $x^*$,  there exists a subgradient $f^{\prime}(x^*)=0$. Then for any $x\in \bbR^d$, one has
\begin{eqnarray}
f(x)&\ge& f(x^*)+\langle  f^{\prime}(x^*), x-x^*\rangle + \frac{\sigma}{s}\|x-x^*\|^s\nonumber\\
&=& f(x^*)+ \frac{\sigma}{s}\|x-x^*\|^s.\label{eq:optim}
\end{eqnarray} 

By Assumption \ref{ass:alg}, one has
\begin{eqnarray}
 f(y_{k+1}) - f(x^*) \le \frac{c_{\cA} \|y_k - x^*\|^{v}}{m_k^{r}}.
\end{eqnarray}
Assume that for some $k_0-2\ge k\ge 0$, \eqref{eq:y-k-R} is valid (it is true for $k = 0$ as $\|x_0-x^*\|_2 \le R$). Then,
\begin{eqnarray}
\frac{\sigma}{s}\|y_{k+1} - x^*\|^s \le f(y_{k+1}) - f(x^*) \le \frac{c_{\cA} R_k^{v}}{m_k^{r}}\le\frac{\sigma}{s2^{s}} R_k^s = \frac{\sigma}{s} R_{k+1}^s
\end{eqnarray}
Thus \eqref{eq:y-k-R} is valid for $k_0-1\ge k\ge 0$. On the other hand,
\begin{eqnarray}
&&f(y_{k+1}) - f(x^*) \le  \frac{c_{\cA}\|y_k-x^*\|^{v}}{m_k^{r}} \le\frac{c_{\cA}\|y_k-x^*\|^{s}R_k^{v-s} }{m_k^{r}}\nonumber\\
&\le&\frac{\sigma}{s2^{s}}\|y_k-x^*\|^s\le\frac{1}{2^s}(f(y_k) - f(x^*)).\label{eq:recur}
\end{eqnarray}
If $k\le k_0$, by recursively using \eqref{eq:recur} from $k-1$ to $1$ and using the fact $f(y_{1}) - f(x^*)\le\frac{\sigma}{s2^{s}}\|y_0-x^*\|^s$ implied by \eqref{eq:recur},  then we obtain \eqref{eq:K-lower}.

\medskip\medskip\medskip\medskip 

\textbf{The case with $k>k_0$.}

\medskip

First by the definition of $k_0$ in Algorithm \ref{alg:restart}, if $k>k_0$, then we have $s<v$. Then by the setting of $k_0$ in Algorithm  \ref{alg:restart} and Theorem \ref{thm:restart}, we have
\begin{eqnarray}
f(y_{k_0}) - f(x^*)\le \frac{1}{2} \left(\frac{\sigma^{v}}{s^v c_{\cA}^s}\right)^{\frac{1}{v-s}}.\label{eq:super-1-v2}
\end{eqnarray}
In \eqref{eq:super-1-v2}, denote $G:=\left(\frac{\sigma^{v}}{s^v c_{\cA}^s}\right)^{\frac{1}{v-s}}.$
For $k\ge k_0$,  by the Step 9 of Algorithm  \ref{alg:restart}, Lemma \ref{lem:simplified-bound} and \eqref{eq:optim}, we have
\begin{eqnarray}
f(y_{k+1}) - f(x^*)&\le& c_{\cA} \|y_k-x^*\|^{v}\nonumber\\
&\le& c_{\cA}\left(\frac{s}{\sigma}(f(y_k) - f(x^*))\right)^{\frac{v}{s}}\nonumber\\
&=&G^{-\frac{v-s}{s}}(f(y_k) - f(x^*))^{\frac{v}{s}}.\nonumber
\end{eqnarray}
Equivalently, we have
\begin{eqnarray}
\frac{f(y_{k+1}) - f(x^*)}{G}\le \left( \frac{f(y_k) - f(x^*)}{G} \right)^{\frac{v}{s}}.\label{eq:super-2-v2}
\end{eqnarray}
Therefore combing \eqref{eq:super-1-v2} and \eqref{eq:super-2-v2}, for $k>k_0$, we have
\begin{eqnarray}
f(y_{k}) - f(x^*)\le \left(\frac{\sigma^{v}}{s^v c_{\cA}^s}\right)^{\frac{1}{v-s}} 2^{-(\frac{v}{s})^{k-k_0}}.
\end{eqnarray}
Therefore, Theorem \ref{thm:restart} is proved.

\medskip\medskip

\end{proof}

\bibliography{references}
\bibliographystyle{siamplain}

\end{document}